\documentclass[11pt, a4paper]{amsart}
\usepackage{a4wide}
\usepackage[pdftex,pdftitle={sub-Riemannian landmark matching as ResNets},bookmarksopen,
colorlinks,
linkcolor=black,
urlcolor=black,
citecolor=black]{hyperref}
\hypersetup{pdfauthor= {E. Jansson and K. Modin}}
\usepackage{makecell}
\usepackage{graphicx,enumitem}
\usepackage{amsmath,amssymb}
\usepackage[english]{babel}
\usepackage{subcaption}
\usepackage{units}
\usepackage{diagbox}
\usepackage{fontenc}
\usepackage{tikz-cd}
\usepackage[numbers]{natbib}
\usepackage{pgfplots}
\usetikzlibrary{calc}
\pgfplotsset{compat=1.16}
\usepackage{xcolor}
\usepackage{caption}
\usepackage{ulem}
\usepackage{cancel}
\usepackage{tikz}
\usepackage{amssymb}

\DeclareSymbolFontAlphabet{\amsmathbb}{AMSb}%
\usepackage{mathbbol}
\usepackage{latexsym}
\usepackage[ruled,vlined]{algorithm2e}
\usepackage{dsfont}

\usepackage{mathtools}

\newcommand{\R}{\amsmathbb{R}}

\newcommand{\N}{\amsmathbb{N}}
\newcommand{\Z}{\amsmathbb{Z}}

\usepackage{booktabs}

\newcommand{\cD}{\mathcal{D}} 

\newcommand{\cF}{\mathcal{F}}

\newcommand{\cH}{\mathcal{H}}

\newcommand{\cL}{\mathcal{L}}
\newcommand{\cM}{\mathcal{M}}

\newcommand{\cR}{\mathcal{R}}
\newcommand{\cS}{\mathcal{S}}

\newcommand{\cU}{\mathcal{U}}

\newcommand{\dd}{\,\mathrm{d}}

\newtheorem{lemma}{Lemma}[section]

\newtheorem{theorem}[lemma]{Theorem}
\theoremstyle{remark}

\theoremstyle{definition}

\usepackage{microtype}

\usepackage{color}
\usepackage[foot]{amsaddr}
\definecolor{darkgreen}{rgb}{0,.6,0}
 
\DeclareMathAlphabet{\mathpzc}{OT1}{pzc}{m}{it} 
\title[Sub-Riemannian Landmark Matching as ResNets]{Sub-Riemannian Landmark Matching and its Interpretation as Residual Neural Networks }

\author[E.~Jansson]{Erik Jansson$^{*,\dagger}$}
\email[]{erikjans@chalmers.se}
	
\author[K.~Modin]{Klas Modin$^{*}$} 
\email[]{klas.modin@chalmers.se}

\address{$^*$Department of Mathematical Sciences, Chalmers University of Technology \& University of Gothenburg, S--412~96 G\"oteborg, Sweden.} 
    
\thanks{
$^\dagger$Corresponding author.
}



\date{\today}

\begin{document}
	\begin{abstract}
    The problem of finding a time-dependent vector field which warps an initial set of points to a target set is common in shape analysis. It is an example of a problem in the diffeomorphic shape matching regime, and can be thought of as a spatial discretization of diffeomorphic image matching. In this paper, we consider landmark matching modified by restricting the set of available vector fields in the sense that vector fields are parametrized by a set of controls.  We determine the geometric setting of the problem, referred to as sub-Riemannian landmark matching, and derive the equations of motion for the controls. We provide two computational algorithms and demonstrate them in numerical examples. In particular, the experiments highlight the importance of the regularization term. A strong motivation is that sub-Riemannian landmark matching have connections with neural networks, in particular the interpretation of residual neural networks as time discretizations of continuous control problems. It allows shape analysis practitioners to think about neural networks in terms of diffeomorphic landmark matching, thereby providing a bridge between the two fields. 
	\end{abstract}

\maketitle
\section{Introduction}
In this paper, we consider the problem of finding an optimal time-dependent vector field transporting a set of initial points on some compact oriented Riemannian manifold $(M,g)$ to some set of targets, i.e.,\ the landmark matching problem of \citet{Joshi2000}.
Landmark matching and, more generally, shape analysis have applications in many image analysis problems. In the field of \emph{computational anatomy} it is used in for instance brain imaging \citep{Bruveris2013,brainC}, MRI scans of the heart \citep{heartmri,heartmri2}, cancer modelling \citep{cancer}, and lung movement studies \citep{lung}.

In diffeomorphic shape matching, shapes are warped by diffeomorphisms generated by smooth vector fields. 
The vector fields are determined by solving an optimization problem. 
Usually, in shape analysis, one optimizes over all vector fields, but in this paper we consider a ``sparse'' setting where only some vector fields are allowed. 
We draw inspiration from the field of neural networks, in which many small components together build up complicated mappings. 
In practice, this means that we consider only a subset of vector fields, parametrized by some set of controls. 
This problem is referred to as \emph{sub-Riemannian landmark matching}~\citep{arguillere2015,Sylvain2016,gris1,Younes2020}.
We furthermore assume that the targets are in a metric space $N$. Points on $M$ are mapped into $N$ by a mapping called the \emph{forward model} (borrowing from the language of inverse problems, \emph{cf.}~\citet{oktem2016}). 
In the paper, we describe the geometry of the sub-Riemannian landmark matching problem and derive the equations of motion for the control parameters. The aim of this paper is to serve as a bridge between shape analysis and the study of neural networks. In particular, we derive from variational principles the equations governing the evolution of the control parameters.

A main point is that sub-Riemannian landmark matching can be thought of as time-continuous neural networks.  
Indeed, residual neural networks (ResNets) can be understood as time discretizations of optimal control problems \citep{Agrachev2021,Celledoni2020,He2015,Vialard2020}. 
Diffeomorphic matching problems are also optimal control problems. This allows us to make connections between shape analysis and neural networks. For instance, just as sub-Riemannian landmarks matching warps points that discretize an image in a low-dimensional space, ResNets   warps high-dimensional input data such as images which, in complete analogy with landmark matching, discretizes a high-dimensional ``meta-image''. 
Furthermore, training the network can be understood as an analogy of the gradient flow method for computing landmark paths.  
The concept of families of vector fields generated by iterated Lie brackets offers an interesting way to view sub-Riemannian landmark matching and, by extension, neural networks.
The richness of the set of mappings that can be obtained is determined by properties of the iterated Lie brackets of the initial vector fields. 

By considering, for instance, so-called \emph{deformation modules}, one can use the sub-Riemannian approach to shape analysis not only as a way to simplify computations, but also as a powerful modelling tool by incorporating knowledge about the desired deformations into the computations \citep{gris1,Younes2020}. 
In this paper, we will approach sub-Riemannian landmark matching from the viewpoint of geometric mechanics, with the goal of being able to describe residual neural networks directly in the presented framework. For this reason, the sub-Riemannian setting is here introduced without making use of deformation modules nor implicit linear constraints on the vector fields.

The inclusion of a forward model has been considered in \citet{oktem2016}. 
Methods from geometry, mechanics, and control have previously been applied to neural networks, both for computational and theoretical aspects \citep{Agrachev2021,nader2021,Ideashape,Tabuada2021UniversalAP,Vialard2020}. 
In this paper, the emphasis is on parameterized landmark matching from the sub-Riemannian perspective.
It is also possible to enhance shape analysis algorithms using neural networks  \citep{DBLP:journals/corr/abs-2102-07951}---a perspective different from ours.  
Further, we mention two papers that pursue goals similar to sub--Riemannian shape analysis.

\citet{Celledoni2023} has utilized deep learning techniques to compose elementary diffeomorphisms into approximations of orientation-preserving diffeomorphisms for use in shape deformation tasks.  

\citet{Scagliotti2023} approximates diffeomorphisms of Euclidean spaces for point cloud matching by taking a control theoretic approach to the composition of linear control systems. 
This view is similar to the examples we present in \autoref{sec:compute_example} and presents an interesting example of how deep learning methods can be applied directly to achieve approximations of diffeomorphisms. 
However, in this paper we pursue a shape perspective, starting in a different geometric setting with the goal to derive equations governing the evolution of the control parameters directly from variational principles. 

A diffeomorphic shape matching perspective in classification problems has been considered in \citet{Younes2020b}, where a kernel-based method inspired by diffeomorphic shape matching 
was used for classification. 
In this paper, we adopt the opposite approach.
We first strive to modify and understand landmark matching from a geometric perspective, and later remark on the connections with neural networks. 
For a general introduction to diffeomorphic shape analysis, see \citet{Younes2010}. 

The paper is structured as follows. 
In \autoref{sec:landmark}, we remind the reader of the geometric approach to shape analysis and landmark matching. 
Sub-Riemannian landmark matching is introduced in \autoref{constrained_landmarks} where we specify the geometric setting  and determine equations of motion for the controls.
A specific parametrization is considered, namely \emph{control-affine systems}. 
In \autoref{sec:compute}, two algorithms for computation of sub-Riemannian landmark matching paths are described. 
Computational examples illustrating some of the theoretical concepts are considered in \autoref{sec:compute_example}.  
In \autoref{sec:resnet}, we describe how shape analysis and neural networks are connected through sub-Riemannian landmark matching.
In \autoref{sec:geometry}, some basic concepts from differential geometry are introduced for the convenience of the reader.
Finally, in \autoref{sec:proof}, an analytical setting for proving the existence of solutions to the sub-Riemannian landmark matching problem is described, after which the existence of solutions is established.

\medskip

\noindent\textbf{Acknowledgements.}
This work was supported by the Wallenberg AI,
Autonomous Systems and Software Program (WASP) funded by the Knut and Alice
Wallenberg Foundation. This work was supported by the Swedish Research Council (grant number 2022-03453) and the Knut and Alice Wallenberg Foundation (grant number WAF2019.0201).
We would like to thank Stefan Sommer for fruitful discussions related to this work.

\section{Diffeomorphic image and landmark matching}\label{sec:landmark}
In this section, we introduce shape analysis and landmark matching. 
Let $(M,g)$ be an orientable compact Riemannian manifold, and let $\operatorname{Diff}(M)$ denote the group of diffeomorphisms of $M$, which is an infinite dimensional Lie group in the Fréchet topology of Sobolev semi-norms (cf.~\cite{Hamilton1982}). 
A diffeomorphism $\varphi \in \operatorname{Diff}(M)$ acts on $f \in C^\infty(M)$ by $\varphi \cdot f=f\circ \varphi^{-1}$.  

Consider now a smooth curve $\gamma(t) \in \operatorname{Diff}(M)$ and denote by $\dot\gamma$ its derivative with respect to~$t$.
A class of weak, right-invariant Riemannian metrics on $\operatorname{Diff}(M)$ is given by
	\begin{align}
	    	\langle \dot{\gamma}, \dot{\gamma} \rangle_\gamma=\int_M v \cdot Lv \, \mathrm{d} x,
    \label{eq:innerp}
	\end{align}
where $v=\dot{\gamma} \circ \gamma^{-1}$ and $L\colon\mathfrak{X}(M) \to \mathfrak{X}(M)$ is an elliptic invertible differential operator.
A typical example is  $L= (1-\alpha \Delta)^k$ for some $k \in \N \cup \{0\} $ and $\alpha>0$. Here, $\Delta$ denotes the Laplace--Beltrami operator acting on vector fields.
	
The \emph{diffeomorphic image matching problem} concerns the matching of two images, $f_0,f_1 \in C^\infty(M)$.  
The matching is performed by minimizing the energy functional over time-dependent vector fields $v\colon [0,1]\to \mathfrak{X}(M)$
	\begin{align}
			\label{eq:enfunc}
			E(v)=\|f_0 \circ \gamma(1)^{-1} -f_1\|^2+\frac{\sigma}{2}  \int_0^1 \int_M v(t) \cdot Lv(t) \mathrm{d} x\mathrm{d}t,
	\end{align}
where $\gamma$ is given by $\dot \gamma(t)=v(t) \circ \gamma(t),\gamma(0)=e$, the identity element of $\operatorname{Diff}(M)$. The parameter $\sigma>0$ is a scalar which determines the regularization strength.

To compute a numerical solution, a  spatial discretization is needed.
One possibility is to use \emph{landmark matching}. The  source  and target images are approximated by a selection of points  $x_1,...,x_m \in M$ that discretizes $f_0$, and $c_1,...,c_m \in M$ that discretizes $f_1$. These points are called landmarks.
The landmark matching problem is to find the time-dependent vector field $v$ solving 
		\begin{align}
		&\!\min_{v} \quad \sum_{i=1}^m d_M^2(y_i(1),c_i)+\frac{\sigma}{2} \int_0^1 \int_M v \cdot Lv \,\mathrm{d}x\mathrm{d}t\label{eq:targetlm}\\
			&\text{such that} \quad \dot y_i(t) = v(t,y_i(t)),~ t\in [0,1], ~y(0)=x_i  \label{eq:conditionlm}.
	\end{align}
Here, $d_M\colon M \times M \to \R$ is the distance function induced by the Riemannian metric on $M$. 

By optimal control theoretic principles, the landmark matching problem can be reformulated in terms of a finite-dimensional Hamiltonian system expressed in the canonical variables $(y_i,p_i)$.
Therefore, landmark matching naturally gives rise to a spatial discretization of the vector fields \citep{Beg2005}. 

Note that a landmark discretization of an image can also be used to track the most important points of interest in the image, for example, the position or outline of facial features in an image of a face. 
In this work, however, we treat landmarks as a way to approximate images. 

\section{Sub-Riemannian landmark matching} \label{constrained_landmarks}
    
In this section, we  consider two modifications of the landmark matching problem described in \autoref{sec:landmark}. 
The first modification is that the vector field  $v$ is parametrized, in the sense that it is determined by some underlying set of control variables $\cU$. 
The available vector fields are restricted to a subset of $\mathfrak{X}(M)$ determined by the parametrization. This subset will \emph{not}  be a Lie subalgebra of $\mathfrak{X}(M)$. 
By working with parametrized vector fields, we are considering a nonlinear control problem where the admissible space of controls are given by a Hilbert space $\cU$.  
We point out that while we present sub-Riemannian landmark matching on a compact manifold $M$, it is possible to work with non-compact manifolds, assuming sufficient decay at infinity of the vector fields. 

The second modification is that the matching term will include the setting when the target landmarks $c_1,\dots, c_m$ are in a metric space $N$ that may be different from $M$. An example is  landmark-label pairs,  in which case $c_1,\dots,c_m$ is in some Euclidean space, the labels being one-hot encoded.  
To achieve this, we supplement the matching term with a function $h: M \to N$ called the \emph{forward model}. The new matching term is 
\begin{align}
    \label{eq:srmatching}
    \sum_{i=1}^m d_N^2(h(y_i(1)),c_i).
\end{align}
Note that if $h$ is the identity map, we retain the original setting.

Diffeomorphic image or landmark matching admits a dynamical formulation, where the source is continuously warped to match the template via a geodesic path of diffeomorphisms.
Constrained landmark matching also admits a dynamical formulation.
Warps here are governed by sub-Riemannian dynamics, thus motivating the designation \emph{sub-Riemannian landmark matching}. 

\begin{figure}[t]
     \centering
     \begin{subfigure}[t]{0.47\textwidth}
        \centering
\includegraphics[width = \linewidth]{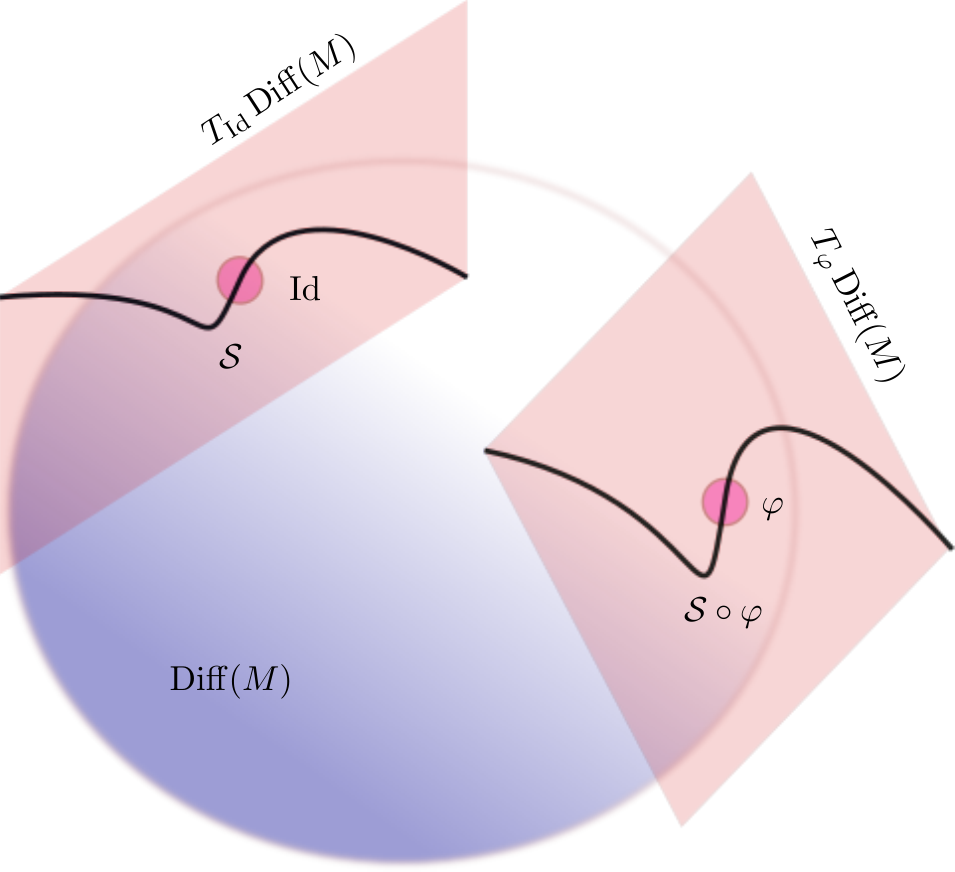}
\caption{General case.}
\label{fig:dist}

 \label{fig:subbundle}
     \end{subfigure}
     ~
     \begin{subfigure}[t]{0.47\textwidth}
         \centering
 \includegraphics[width =\linewidth]{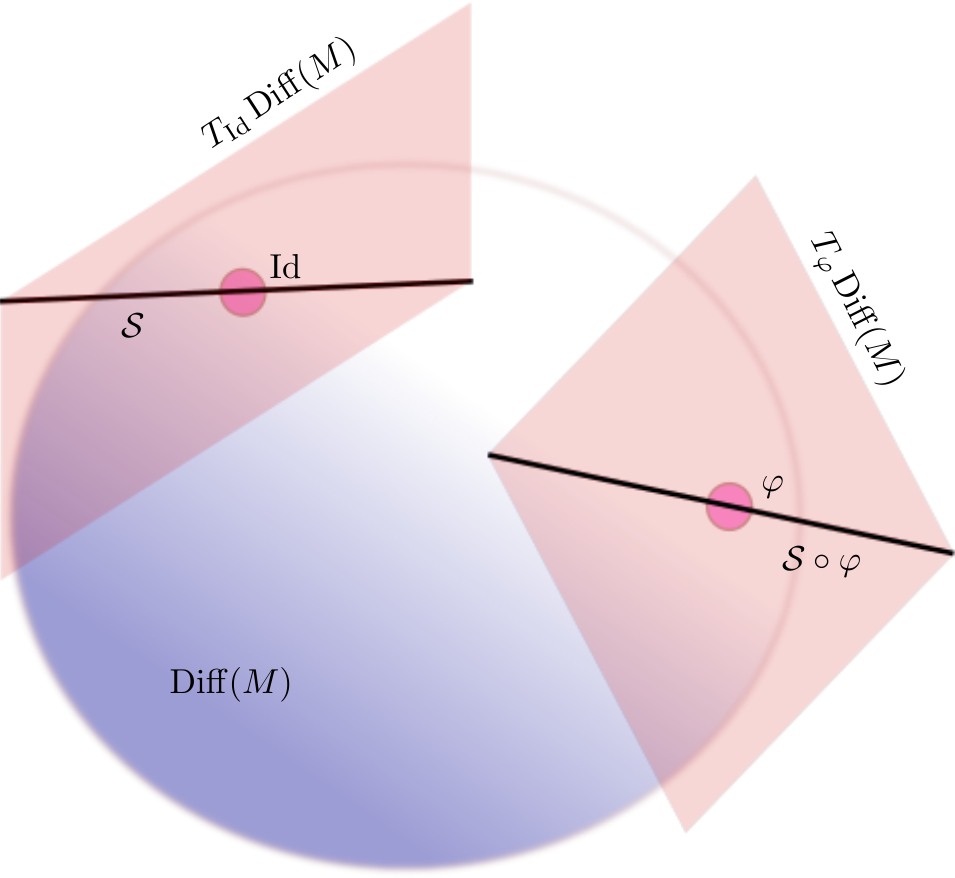}
\caption{Control-affine case.} 
\label{fig:subbundle_dist}
     \end{subfigure}
    \caption{Illustration of 
$\mathcal E(\mathcal S)$  sub-bundle of $T\operatorname{Diff}(M)$ in two different cases. 
}
    \label{fig:bundle}
\end{figure}

Consider a submanifold of vector fields $\mathcal{S} \subset\mathfrak{X}(M)$.
By right translation, $\mathcal{S}$ gives rise to a right-invariant sub-bundle of $T\operatorname{Diff}(M)$ in the category of fiber bundles, see \autoref{fig:subbundle}. We denote it $\mathcal E(\mathcal S)$. 
The fiber above $\varphi\in\operatorname{Diff}(M)$ is the submanifold of $T_\varphi\operatorname{Diff}(M)$ given by $\mathcal{S}\circ\varphi$. 
The submanifold $\mathcal S$ is parameterized by some function $F$ mapping from parameters into vector fields,
\begin{equation}
        v \in \mathcal S \iff
    \exists u \in \cU \text{ s.t. } v(x) = F(u)(x). \label{eq:S}
\end{equation}

 The minimization problem is therefore 
\begin{align}
    	&\!\min_{u\colon[0,1] \to\, \cU } \quad \sum_{i=1}^m d_N^2(h(y_i(1),c_i))+\int_0^1 \ell(F(u)) \mathrm{d}t\label{eq:targetvlm_nonformal}\\
			&\text{s.t.  } \quad \dot y_i =  F(u)(y_i), ~t\in [0,1], ~y_i(0)=x_i \label{eq:constvlm_nonformal}.
\end{align}
where the regularization term is determined by the function $\ell: \mathfrak{X}(M) \to \R$.   

We continue by studying the dynamics of the sub-Riemannian landmark matching problem. \autoref{eq:targetvlm_nonformal} is reminiscent of an action functional in Lagrangian mechanics. 
In fact, since the shape matching term, \autoref{eq:srmatching}, depends solely on the final point of the landmark paths, an optimal $u$ must follow the dynamics determined by the action functional consisting only of the regularization term.
We make the additional assumption that 
\begin{align*}
    \ell(v)=\frac{\sigma}{2} \langle v,v \rangle 
\end{align*}
where $\langle \cdot,\cdot\rangle$ is the inner product defined in \autoref{eq:innerp}. 
The dynamics is unaffected by $\sigma$, so in the derivations below we take $\sigma =1$ without loss of generality.

The term $\ell(v)$ is the \emph{reduced Lagrangian} corresponding to the full Lagrangian on $T\operatorname{Diff}(M)$ given by $L(\gamma,\dot\gamma) = \ell(\dot\gamma\circ\gamma^{-1})$.
With this observation in mind, we can proceed as in the proof in \citet[Theorem 13.5.3]{Marsden1999}. 
First, for a variation $\gamma_\epsilon(t)$, we have
\[
    \left.\frac{\mathrm d}{\mathrm d\epsilon}\right|_{\epsilon = 0} \dot\gamma_\epsilon\circ\gamma_\epsilon^{-1} = \dot w - [v, w]
\]
where $v = \dot\gamma\circ\gamma^{-1}$ and $w = \mathrm d\gamma_{\epsilon}/\mathrm d\epsilon|_{\epsilon = 0}   \circ\gamma^{-1}$ such that $w \in T_v \mathcal{S}$.
Noting that 
$\frac{\delta \ell}{\delta v}=v$ the variational derivative of the energy functional, $\delta E$, is given by 
\begin{align}
    \label{eq:deltaE}
     \delta E= \int_0^1 \left\langle v , \delta_u v \right\rangle \mathrm{d}t =  \int_0^1 \left\langle v , \dot w -\text{ad}_v w \right\rangle \mathrm{d}t=\int_0^1 \left\langle -\dot v- \text{ad}_v^T v  , w \right\rangle \mathrm{d}t,
\end{align}
where the convention that  $\text{ad}_v w = [v,w]$ is used. Further, $\text{ad}_v^T$ denotes the transpose of $\text{ad}_v$ with respect to the inner product. 
By calculus of variations, the condition $\delta E =0$ is equivalent to
\begin{align}
    \label{eq:eup}
    \dot{v} = \text{ad}_v^T v +\cM,
\end{align}
where $\cM \in T^\perp_{v}\cS = \{w_1 \in \mathfrak{X}(M) | \langle w_1,w_2 \rangle = 0 \text{ for all } w_2 \in T_v \cS \}$ can be thought of as a Lagrange multiplier to ensure that we remain on $\cS$.

We now assume that the vector fields take values in a Sobolev space $\mathsf{H}^{s}(M)$ and that $L\colon\mathsf H^{s}(M) \to \mathsf H^{s-2k}(M)$. 
Here, $s>\frac{n}{2}+1$ so that $\mathsf{H}^s(M)$ is continuously embedded into the space of continuously differentiable vector fields, \citep[Section 2.2]{Bruveris} meaning that for any $v \in \mathsf H^{s}(M)$ there is a $C>0$ such  that 
\begin{align*}
    \|v\|_{1,\infty} \leq C\|v\|_{\mathsf H^{s}},
\end{align*}
where 
\begin{align*}
    \|v\|_{1,\infty}=\sup_{x \in M }\left( |v(x)|+\sum_{i=1}^{n} |\nabla v^i(x)|\right).
\end{align*} 
The space of continuously differentiable vector fields is denoted by  $C^1(M;TM)$. 

With $m = L v$, \autoref{eq:eup} is equivalent to    
\begin{align}
    \label{eq:eup-proj}
    \dot{m} = \nabla_v^T m +\nabla_m v +m\operatorname{div}(v) + L \cM,
\end{align}
where $\nabla$ denotes the covariant derivative and $\nabla^T$ its $\mathsf{L}^2$ adjoint \citep{Bruveris2013}.
Let $D_um$ denote the derivative of $m$ with respect to $u$.
By the chain rule, we obtain
 \begin{align}
     \dot{m} = D_u m \cdot \dot{u},
 \end{align}
and so \autoref{eq:eup-proj} becomes
\begin{align}
    \label{eq:eup_onestep}
    D_u m \cdot \dot u  = \nabla_v^T m +\nabla_m v +m\operatorname{div}(v) + L \cM.
\end{align}
Thus, with $A(u) \coloneqq D_u^T m \circ L^{-1} \circ D_u m  $ the equation of motion of $u$ is
\begin{align}
    \label{eq:eup-proj_newandimproved}
    A(u)\dot u = (D_u^T m \circ L^{-1}) \left(\nabla_v^T m +\nabla_m v +m\operatorname{div}(v)\right),
\end{align}
since $\cM \in  T^\perp_{v}\cS$.

As a concrete example, let $\cU= \R^l$ and assume that the vector fields are given by 
\begin{align}
    \label{eq:linaff}
    F(u)=X^0 + \sum_{i=1}^l u_i X^i,
\end{align}
where $X^i \in \mathfrak{X}(M)$.
Vector fields of the form~\eqref{eq:linaff}, known as control--affine systems, are often studied in geometric control theory \citep{Jurdjevic1999}. 

In this setting, the geometry is simplified.
In particular,  $\mathcal{E}(\cS)$ of \autoref{constrained_landmarks} becomes an \emph{affine distribution}. 
Note that if $X^0 =0$, this will be a linear sub-bundle of $T\operatorname{Diff}(M)$. This is due to the linear-affine structure of  \autoref{eq:linaff} \citep{Clelland2009}. See  \autoref{fig:subbundle_dist} for an illustration.

We continue with the dynamics of  $u$  in the setting of \autoref{eq:linaff}, with $X^0= 0$. 
This amounts to calculating the projection $D_u^T m$.
First, we note that \begin{align*}
    D_u m : \R^l \to \mathfrak{X}(M)^*, 
\end{align*}
and that $D_u m \cdot \dot{u} =\sum_{i=1}^l L X^i \dot{u}_i$.
Thus, for an arbitrary $w \in \mathfrak{X}(M)$, 
\begin{align*}
    &\langle D_u  m \cdot \dot{u}, w \rangle_{ \mathfrak{X}(M)^*, \mathfrak{X}(M)} \\
   &~~=  \left\langle \sum_{i=1}^l LX^i \dot u_i, w \right\rangle_{\mathsf{L}^2(M)}= \sum_{i=1}^l \dot u_i \int_{M} L X^i \cdot w \dd x= \langle \dot{u},D_u^T m \cdot w \rangle_{\R^l ,\R^l},
\end{align*}
meaning that  $D_u^T m: \mathfrak{X}(M) \to \R^l$ evaluated at $w$  is given by 
\begin{align*}
   D_u^T m \cdot w  =\sum_{i=1}^l \int_{M} LX^i \cdot w \dd x \; e_i ,\end{align*}
where $(e_i)_{i=1}^m$ is the standard $\R^l$ basis.  
As $A(u) =D_u^T\circ L^{-1} \circ D_u m$,
\begin{align*}
    A(u) \dot{u} =\sum_{i=1}^l \sum_{j=1}^l \dot{u}_j  \int_{M} LX^i \cdot X^j  \dd x \; e_i,
\end{align*}
so the operator $A(u):\R^l \to \R^l $ is an $m\times m$ matrix with entries
\begin{align*}
    [A(u)]_{ij} = \int_{M} LX^i \cdot X^j.
\end{align*} 
Due to the linear structure of $\mathcal{E}(\cS)$, $A$ is in this case independent of $u$.

A natural question is how  the vector fields $X^0, \dots,X^l$ should be chosen. One possible approach is to set $X^0=0$ as we did above and choose the remaining vector fields such that the Lie algebra of vector fields generated by the fields is made as large as possible.
Concepts from the field of geometric nonlinear control offer a way to make this precise, and is an interesting way to view sub-Riemannian landmark matching. 
As an example, given a family of smooth vector fields $\cF$, consider  
$\operatorname{Lie}(\cF)$,
the smallest Lie subalgebra containing $\cF$. One can show that this will be spanned by the iterated Lie brackets of the vector fields in $\cF$ and we say that the family is \emph{bracket generating} for $\operatorname{Lie}(\cF)$ \citep{Agrachev2004,Hrmander1967}. 
In control-affine systems, the family of vector fields described by \autoref{eq:linaff} is bracket generating if  $\{X^1,...,X^l\}$ is.  The notion of integrability of the distribution generated by $\cF$ gives an intuition as to how controllability can be interpreted and why sub-Riemannian landmark matching works. If the distribution is integrable, then at every point it will be tangent to the leaves of a foliation, meaning that if we only can select vector fields from $\cF$ we are stuck in the leaf of the initial point $x_0$; we cannot connect any two points by an integral curve unless they are in the same leaf. 
The condition that the family is bracket generating tells us that the distribution behaves in a sense \emph{oppositely} to this. 
A guiding principle for selecting the vector fields is therefore to make sure that $\operatorname{Lie}(\cF)$ is as large as possible.  
Further, we remark the solution of differential equations such as \autoref{eq:constvlm_nonformal} with $F(u)$ given by \autoref{eq:linaff} with $X^0=0$ can be formally described by series containing the iterated Lie brackets of the vector fields $X^1,\dots, X^l$, so another important goal when choosing the vector fields is that their iterated Lie brackets should generate a large class of vector fields. 
\citep{Beauchard2023,Sussmann1986}

\section{Computation of sub-Riemannian landmark matching } \label{sec:compute}
In the case of diffeomorphic landmark matching, various computational algorithms such as shooting or gradient flow-based methods are available. 
Both of these approaches are also viable for computing sub-Riemannian landmark paths.  

We first note that the energy is conserved along the dynamics. 
Indeed, it holds that 
\begin{align*}
    \frac{\mathrm d \ell(v)}{\mathrm d t} = \langle v, \dot v \rangle =  \langle v, \operatorname{ad}_v^T v + \cM \rangle = \langle \operatorname{ad}_v^v,v\rangle = 0.
\end{align*}
Therefore, the energy functional can be written
\begin{align*}
    \tilde{E}=\frac{\sigma}{2} \langle F(u_0),F(u_0)\rangle  +\sum_{i=1}^l d_N^2(h(y_i(1)),c_i).
\end{align*}
In the case of control-affine systems with $X^0 = 0 $, the regularization term is given by 
\begin{align*}
    \frac{\sigma}{2} \langle F(u_0),F(u_0)\rangle = \frac{\sigma}{2}\sum_{i,j=1}^l u_i u_j \int_M LX^i \cdot X^j \dd x,
\end{align*}
and is particularly simple if the vector fields are orthogonal eigenfunctions of $L$.  

The dynamics of the control variable, if known, taken together with the landmark dynamics, constitute a dynamical system which can be numerically integrated to approximate $y_i(1)$ so that one may evaluate $\tilde{E}$.  

Then, we optimize the initial value of the control, $u_0$, using for instance gradient descent, 
\begin{align}
    \label{eq:shooting}
    u_0^{[j+1]} = u_0^{[j]}-\epsilon \nabla_{u_0^{[j]}} \tilde{E},
\end{align}
where $\epsilon$ denotes the gradient descent step size and the gradient of 
$\tilde{E}$ with respect to $u_0$ can be computed using automatic differentiation. 
One can use, for instance, PyTorch in Python or ForwardDiff package in Julia \citep{NEURIPS2019_9015,RevelsLubinPapamarkou2016}.  
However, we are not limited to  gradient descent. 
Indeed, any method for minimizing $\tilde{E}$, considered a function of $u_0$ can be used as long as the computation of the gradient is feasible. 
In practice, packages that combine optimization routines with the computation of gradients (or higher-order derivatives) using automatic differentiation can be used, for instance the Julia package Optim~\citep{Optim.jl-2018}.

In cases where the dynamics of the control variables are  difficult to compute explicitly, we can take a different approach and work with a discretization of the entire control path,  $u_{t_0},...,u_{t_n}$, in a partition of the time interval, $0=t_0 <t_1 <...<t_n =1$. 
One may then view $\tilde{E}$ as a function of $u_{t_0},...,u_{t_n}$
and optimize each of these variables using, for instance, gradient descent. 
\begin{align}
    \label{eq:graddiff}
    u_{t_k}^{[j+1]} = u_{t_k}^{[j]}-\epsilon \nabla_{u_{t_k}^{[j]}} \tilde{E}. 
\end{align}
Note that we must perform one additional gradient descent step for each point in the discretization of the time interval. Knowing the dynamics of $u_t$ therefore reduces the number of variables to one. 
This would seemingly imply that a shooting method is always preferable, but we emphasize that in numerous instances an explicit ODE for $u$ cannot be written down, for instance, when $A(u)$ cannot be given explicitly.

\section{Numerical examples}
\label{sec:compute_example}
In this section, we consider several concrete numerical examples of sub-Riemannian landmark matching on the flat torus $\R^2 / (2\pi\Z)^2$.

\subsection{Landmark deformations without forward model}
\begin{figure}[ht]
     \centering
     \includegraphics[width=0.6\textwidth]{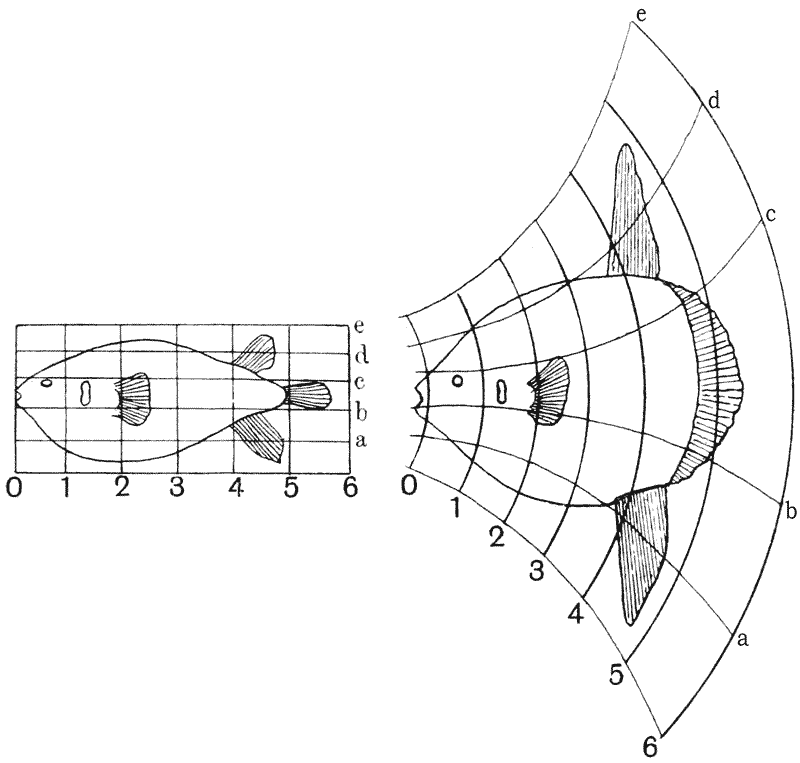}
    \caption{The left fish is the initial image, and the right one is the target. Note the structural similarities between the fish, implying that it should be possible to warp one to the other. }
    \label{fig:fish}
\end{figure}

\begin{figure}[!htb]
     \centering
     \begin{subfigure}[t]{0.4\textwidth}
         \centering
         \includegraphics[width=\textwidth]{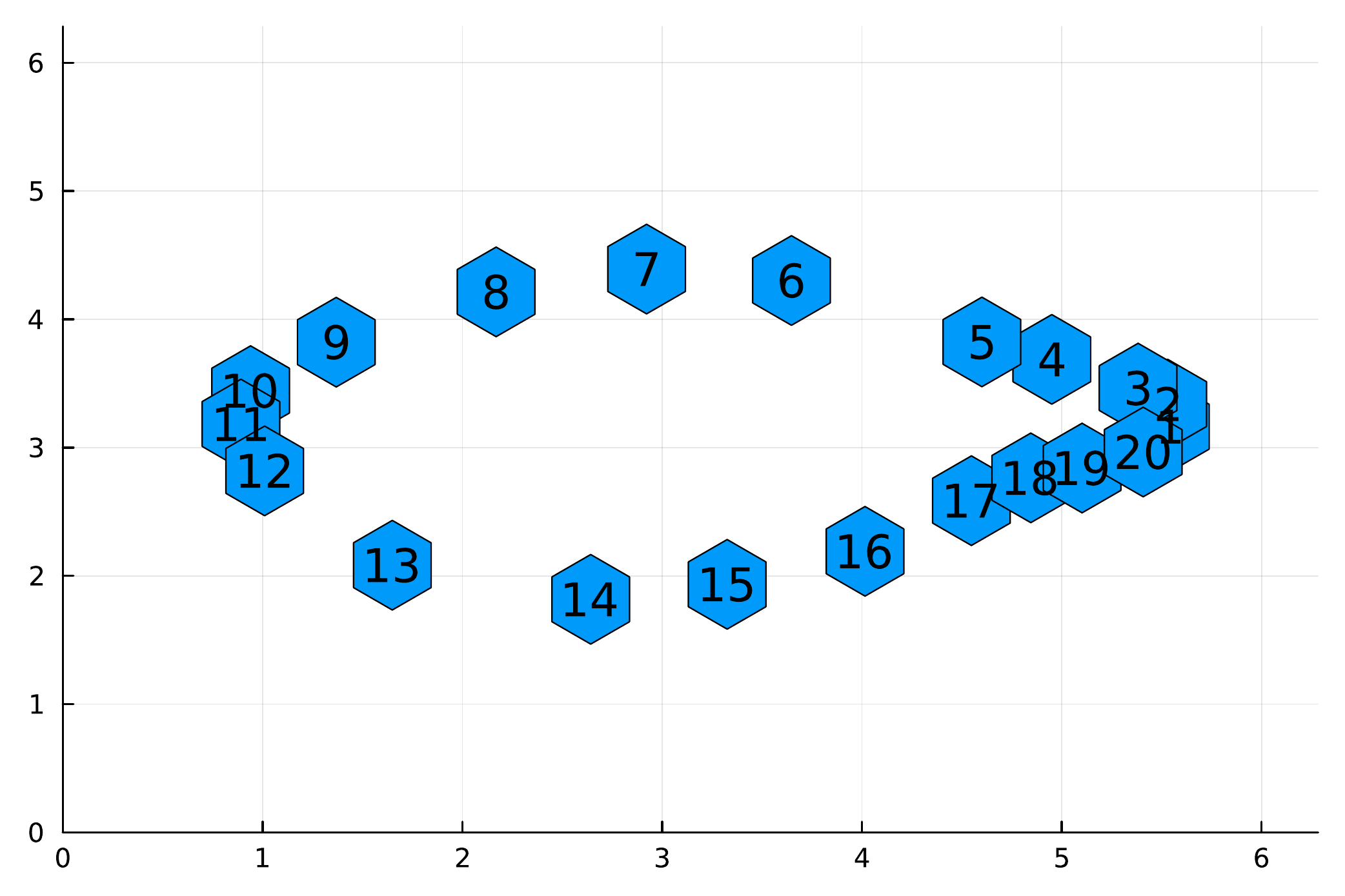}
         \caption{Initial landmarks.}
         \label{fig:fish1_lm}
     \end{subfigure}
     ~
     \begin{subfigure}[t]{0.4\textwidth}
         \centering
         \includegraphics[width=\textwidth]{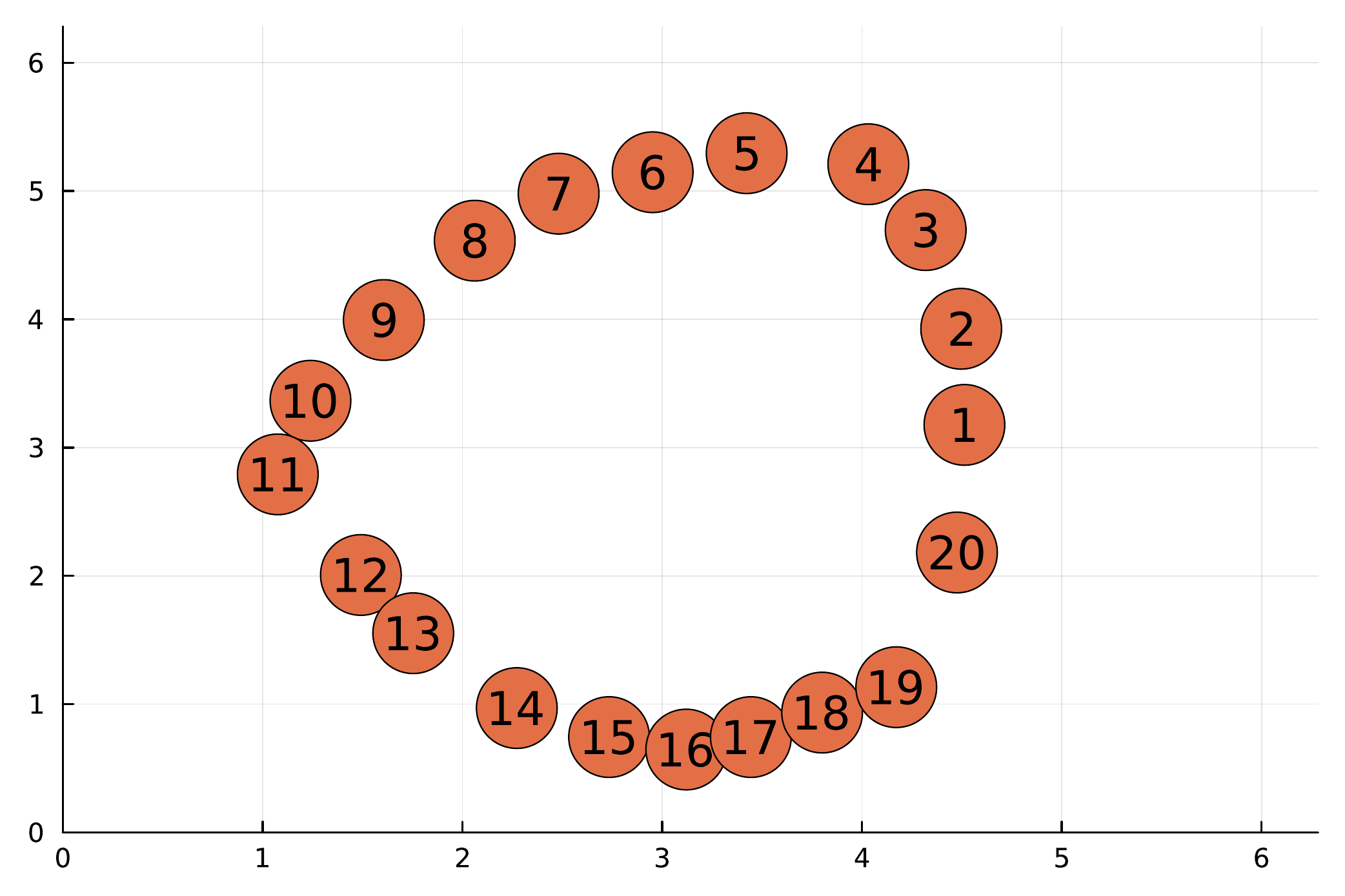}
         \caption{Target landmarks.}
         \label{fig:fish2_lm}
     \end{subfigure}
     \label{fig:fish_lm}
     \caption{The images in \autoref{fig:fish} approximated by landmarks on the flat torus.}
\end{figure}

We first perform the matching on the classical images of fish from D'Arcy Thompson's classic book~\citep{Thompson1992} (see \autoref{fig:fish}). 
The two images in \autoref{fig:fish} are discretized using labelled landmarks. 
The goal is to warp the landmarks approximating the first image, see \autoref{fig:fish1_lm}, to those approximating the second image, see \autoref{fig:fish2_lm} in the setting of  \autoref{eq:linaff} with various choices of the dimension $l$ of $\cU = \R^l$ and the set of vector fields $X^1, \dots, X^l$. 
Thereby, we can illustrate the discussion at the end of \autoref{sec:compute} on how the behavior of the iterated brackets allow for more complicated warps. 

\subsubsection{Trivial control dynamics and trivial landmark translations}
We begin with $l=2$.  Let us start by choosing $X^1 = e_1$ and $X_2 = e_2$, where $e_1$ and $e_2$ are the standard basis vectors. These vector fields are orthogonal, so $A(u)$ is invertible. 
Note that the dynamics of $u = (u_1,u_2) \in \cU$ are trivial, i.e., $\dot u = 0$ implying that $u(t) = u(0)$ for all times. 
Further, the Lie bracket of these vector fields is zero, so the best mapping we can hope for is a linear combination of $e_1$ and $e_2$, or, in other words, a translation.

We use the shooting method outlined in \autoref{sec:compute} to find the optimal $u_0$. 
We use Julia with the package DiffEq \citep{rackauckas2017differentialequations} to approximate solutions to the ordinary differential equations. 
Further, we perform the optimization of $u_0$ with the package Optim. 
We use the L-BFGS method, with gradients computed using ForwardDiff.
The regularization strength is set to $0$ and an integration step size of $h = 10^{-1}$ is used.
The initial guess of $u_0$ is the zero vector.

\begin{figure}[!htb]
     \centering
     \begin{subfigure}[t]{0.48\textwidth}
         \centering
         \includegraphics[width=\textwidth]{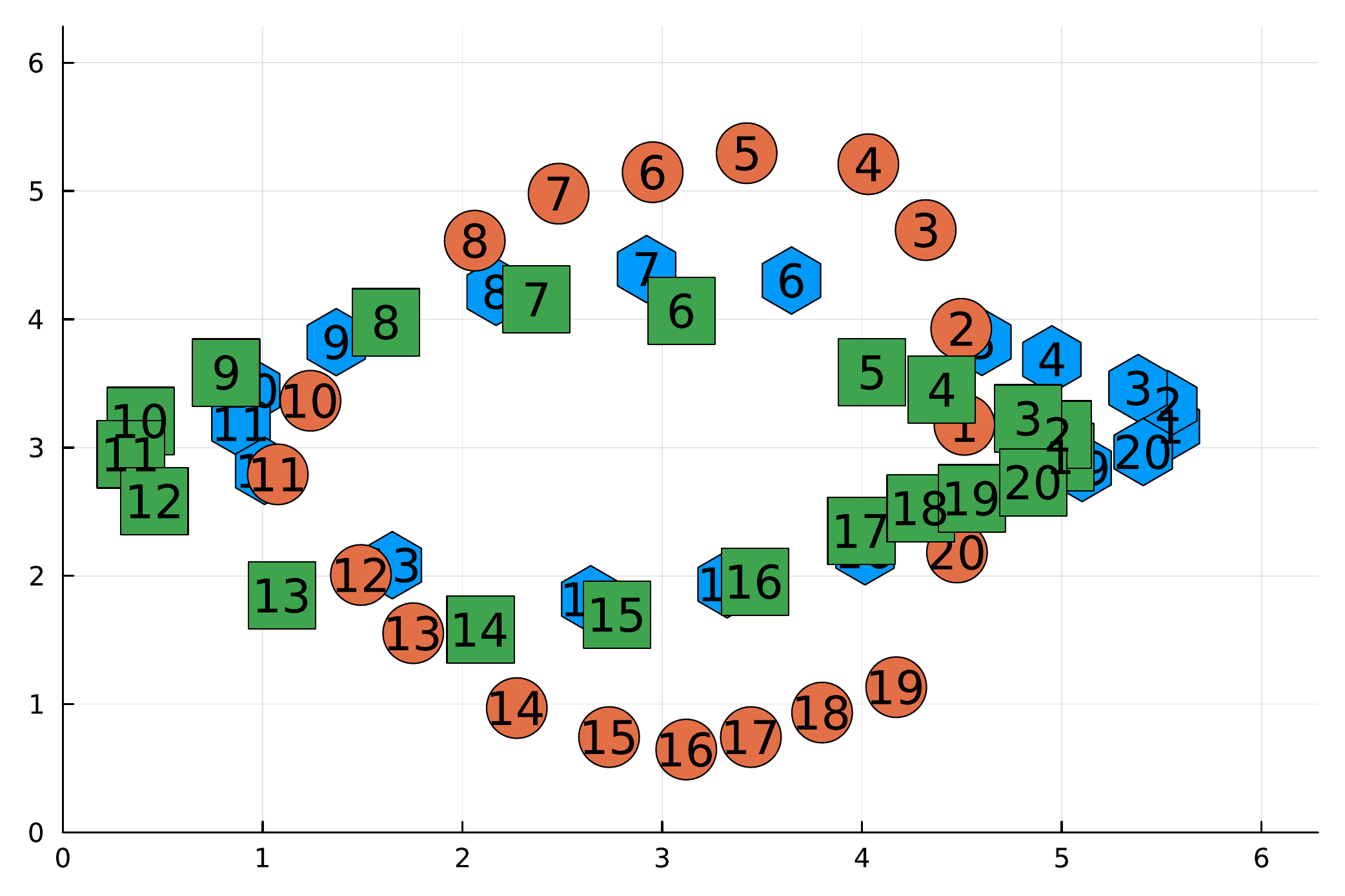}
         \caption{The blue hexagons are the initial landmarks, the orange circles are the targets and the green squares are the deformed landmarks. }
         \label{fig:result_trivial}
         \end{subfigure}
         ~
         \begin{subfigure}[t]{0.48\textwidth}
         \centering
           \includegraphics[width=\textwidth]{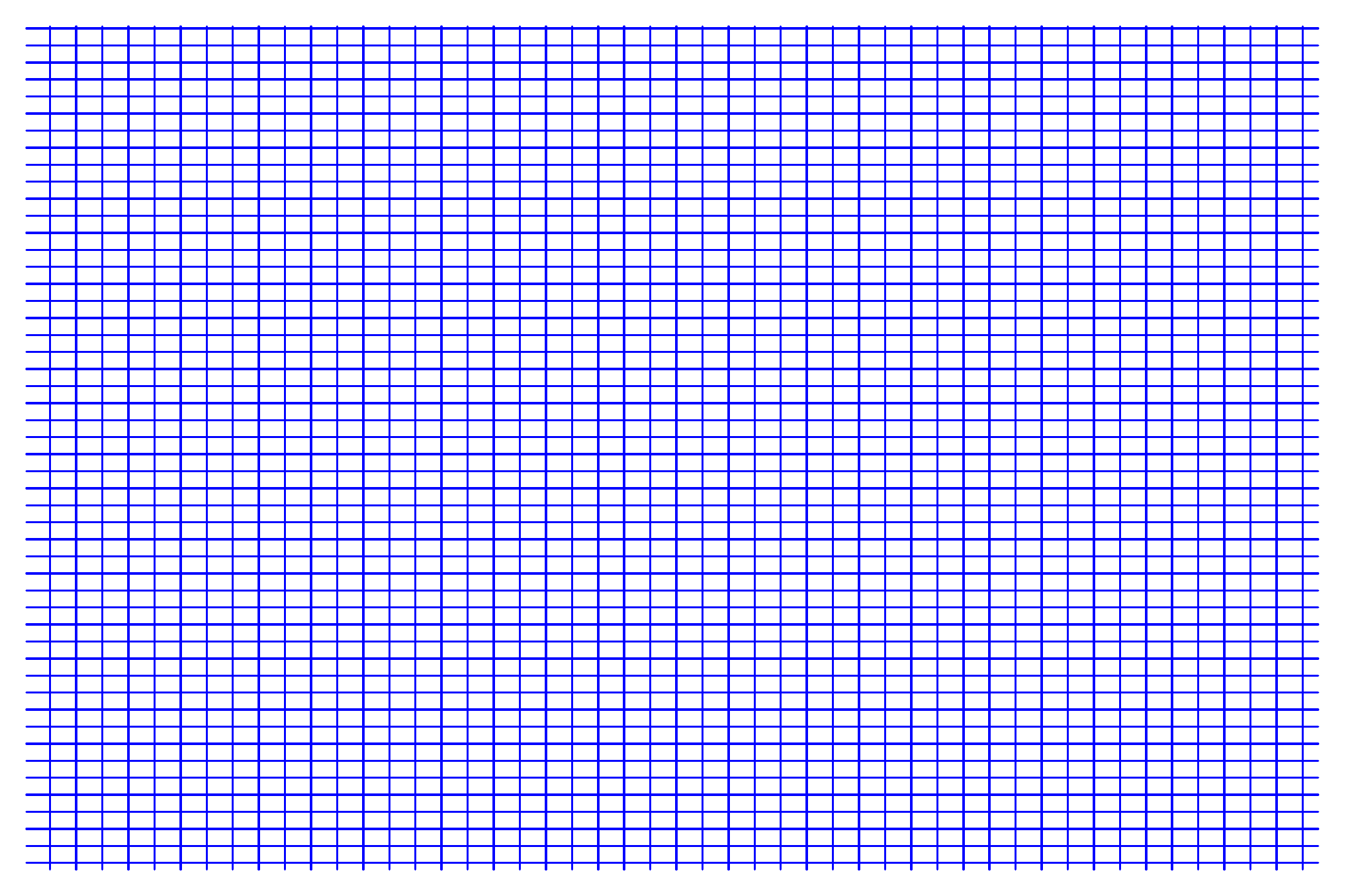}
         \caption{ Grid deformation illustrating the warp. }
         \label{fig:warp_trivial}
         \end{subfigure}
         \caption{The result of applying the shooting method with constant vector fields. We illustrate both how the initial landmarks are transported and the resulting grid deformation. The deformed landmarks are only a translation of the original landmarks.}
         \label{fig:trivial}
\end{figure}

The optimization yields the moved landmarks in  \autoref{fig:result_trivial}. In \autoref{fig:trivial} the warp  is visualized as a grid deformation. 
As expected, we only obtain a translation of the landmarks, as the selected vector fields are unable to produce more advanced warps. 

\subsubsection{Non-trivial control dynamics and almost-trivial landmark translations}
In an attempt to obtain more complicated warps, we instead select $X^1 = e_1$, $X^2 =  xe_1$, $X^3 = e_2$ and $X^4 = xe_4$ so that $\cU = \R^4$.

\begin{figure}[!htb]
     \centering
     \begin{subfigure}[t]{0.48\textwidth}
         \centering
         \includegraphics[width=\textwidth]{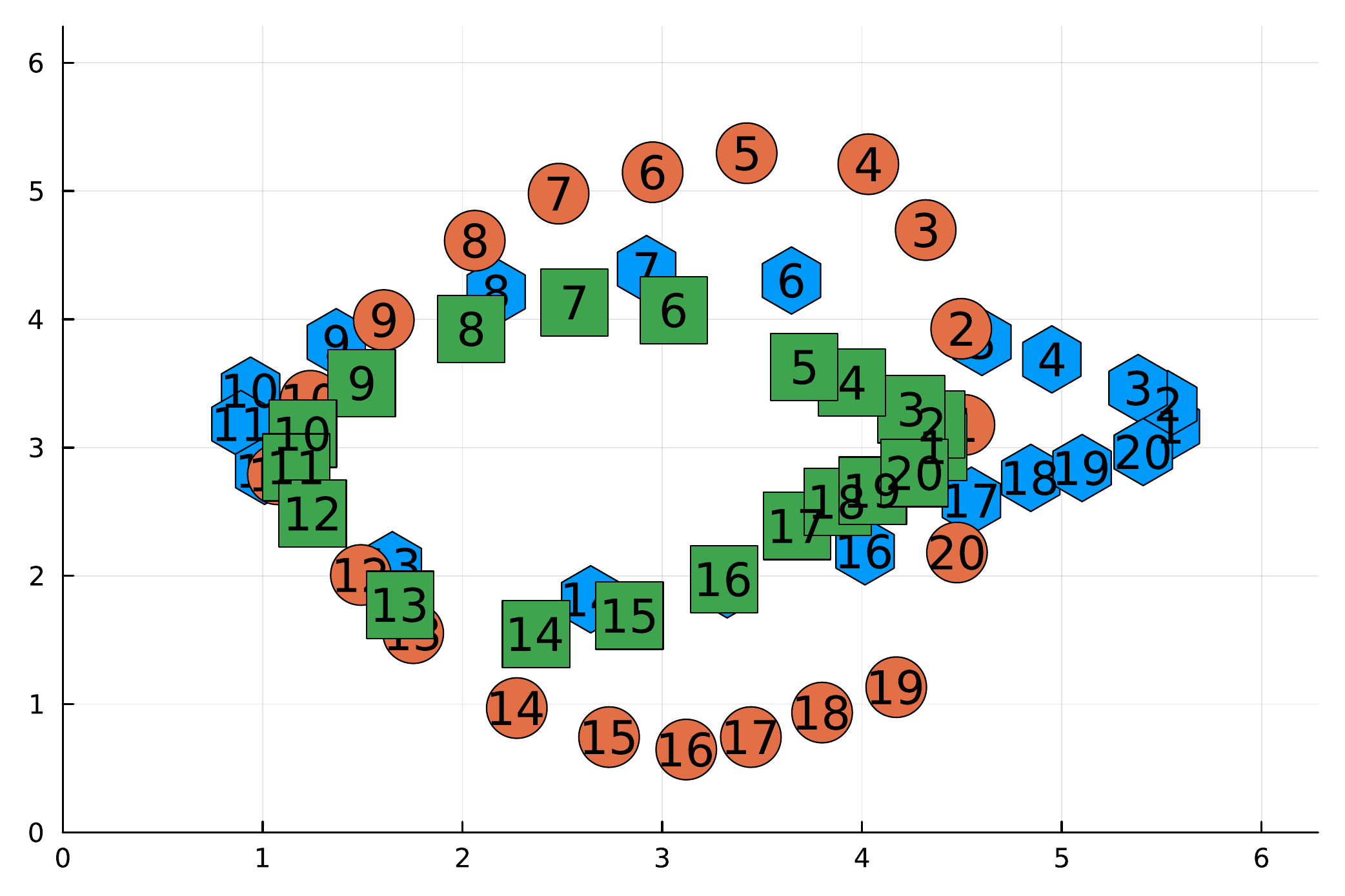}
         \caption{The blue hexagons are the initial landmarks, the orange circles are the targets and the green squares are the deformed landmarks. }
         \label{fig:result_almost_trivial}
         \end{subfigure}
         ~
         \begin{subfigure}[t]{0.48\textwidth}
         \centering
           \includegraphics[width=\textwidth]{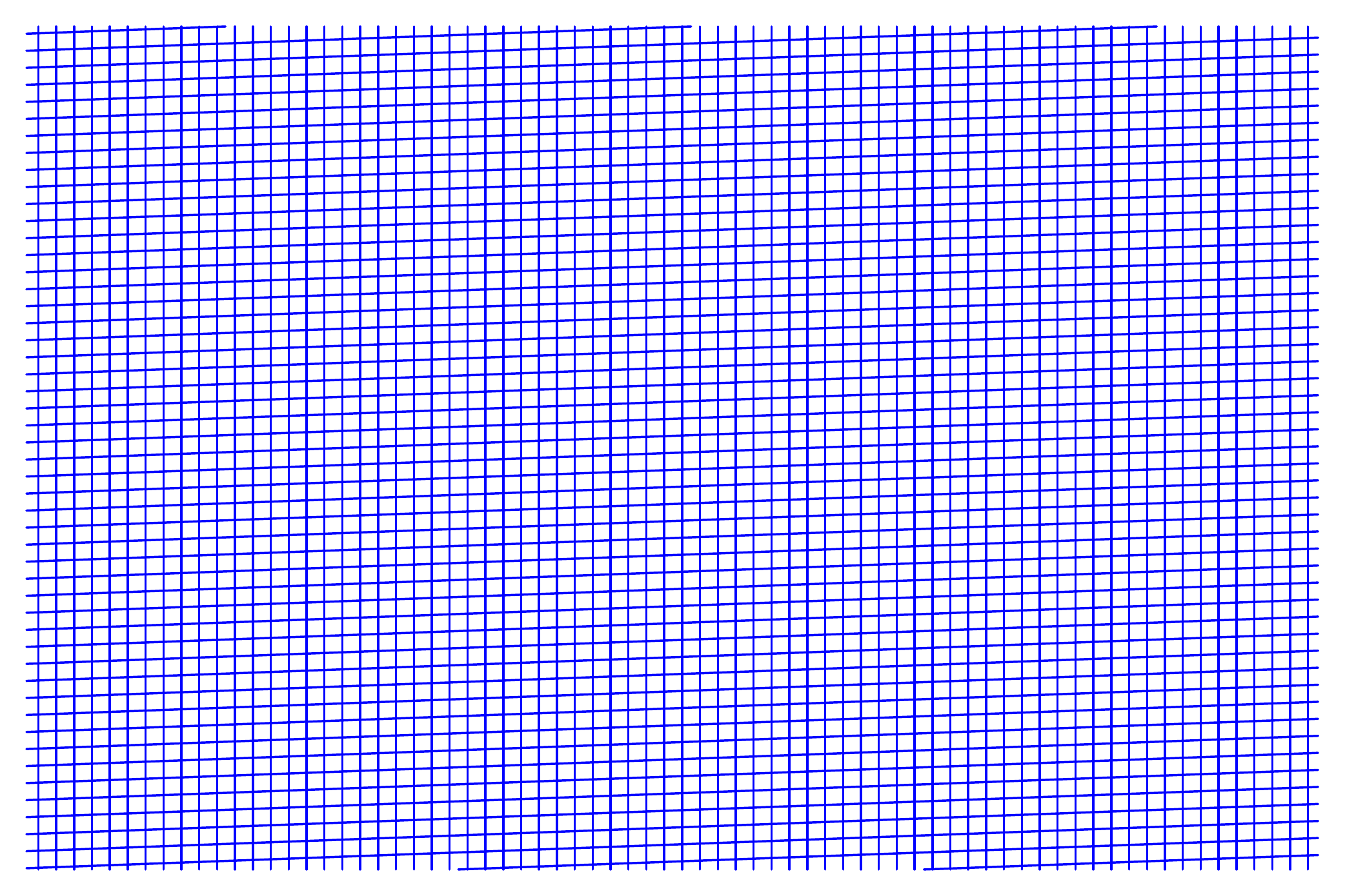}
         \caption{ Grid deformation illustrating the warp. }
         \label{fig:warp_almost_trivial}
         \end{subfigure}
         \caption{The result of applying the shooting method with  four vector fields that form an algebra. }
         \label{fig:almost_trivial}
\end{figure}

The matrix $A(u)$ is invertible, and we rewrite  \autoref{eq:eup-proj_newandimproved} as
\begin{align}
    \label{eq:tordyn} 
    \dot{u} =  A(u)^{-1}(D_u^T m \circ L^{-1}) \left(\nabla_{v}^T m + m \operatorname{div}(v) +\nabla_m v  \right). 
\end{align}
By writing \autoref{eq:tordyn} in coordinates one obtains the equations  for the evolution of $u$,

\begin{align*}
&\dot{u}_1 = 3u_1u_2 + u_3u_4, \quad \dot{u}_2 = 3u_2^2 + u_4^2, \\
&\dot{u}_3 = u_2u_3 + u_1u_4, \quad ~\, \dot{u}_4 = 2u_2u_4. 
\end{align*}
With all other parameters and methods as in the previous example, we obtain the moved landmarks in \autoref{fig:almost_trivial}.
In \autoref{fig:warp_almost_trivial} the warp  is visualized as a grid deformation. 

 Note that the deformed landmarks are obtained by translating  the initial landmarks and slightly compressing them in the $x$-direction, so even though we have non-trivial dynamics of $u$, we do not obtain a complicated warp. 
 This can be understood by considering the Lie brackets of the vector fields. 

Indeed, it holds that the set $\left\{X^1,X^2,X^3,X^4\right\}$ forms a Lie subalgebra of $\mathfrak{X}(\R^2/(2\pi \Z)^2)$, and in the light of the discussion at the end of \autoref{sec:compute}, the set $\left\{X^1,X^2,X^3,X^4\right\}$ is a poor choice of vector fields because we do not obtain a rich family of iterated Lie brackets that can produce complicated warps. 
 
 Further, since the vector fields form an algebra, including more steps in the integration scheme is not going to result in more complicated warps. 
Indeed, reducing, or increasing, the integration step size does not alter the outcome depicted in \autoref{fig:almost_trivial}.

\subsubsection{Non-trivial control dynamics and non-trivial landmark translations}
To obtain a more interesting warp, we select the vector fields such that they do not form a Lie subalgebra. 
To this end, we pick $X^1 = \sin(x)e_1$, $X^2 = \sin(y/2)e_1$, $X^3 = \sin(x/2)e_2$ and $ X^4 = \sin(y/4) e_2$. 
Here, the iterated Lie brackets of $X^1,X^2,X^3$ and $X^4$ are non-trivial. Thus, we should be able to obtain more elaborate warps than in the previous example. 

\begin{figure}[!htb]
     \centering
     \begin{subfigure}[t]{0.48\textwidth}
         \centering
         \includegraphics[width=\textwidth]{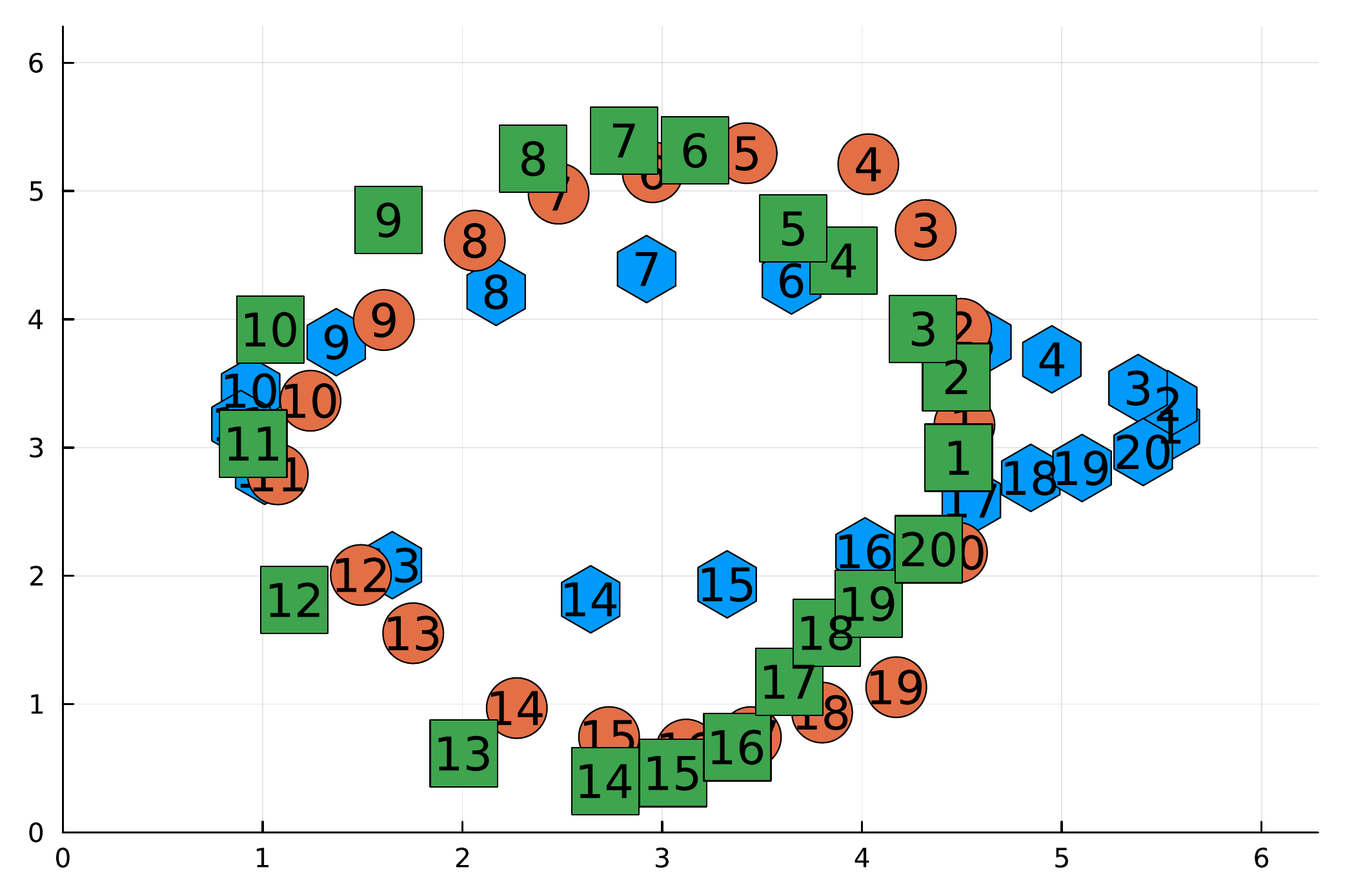}
         \caption{The blue hexagons are the initial landmarks, the orange circles are the targets and the green squares are the deformed landmarks.}
         \label{fig:4d_h01.pdf}
         \end{subfigure}
         ~
         \begin{subfigure}[t]{0.48\textwidth}
         \centering
           \includegraphics[width=\textwidth]{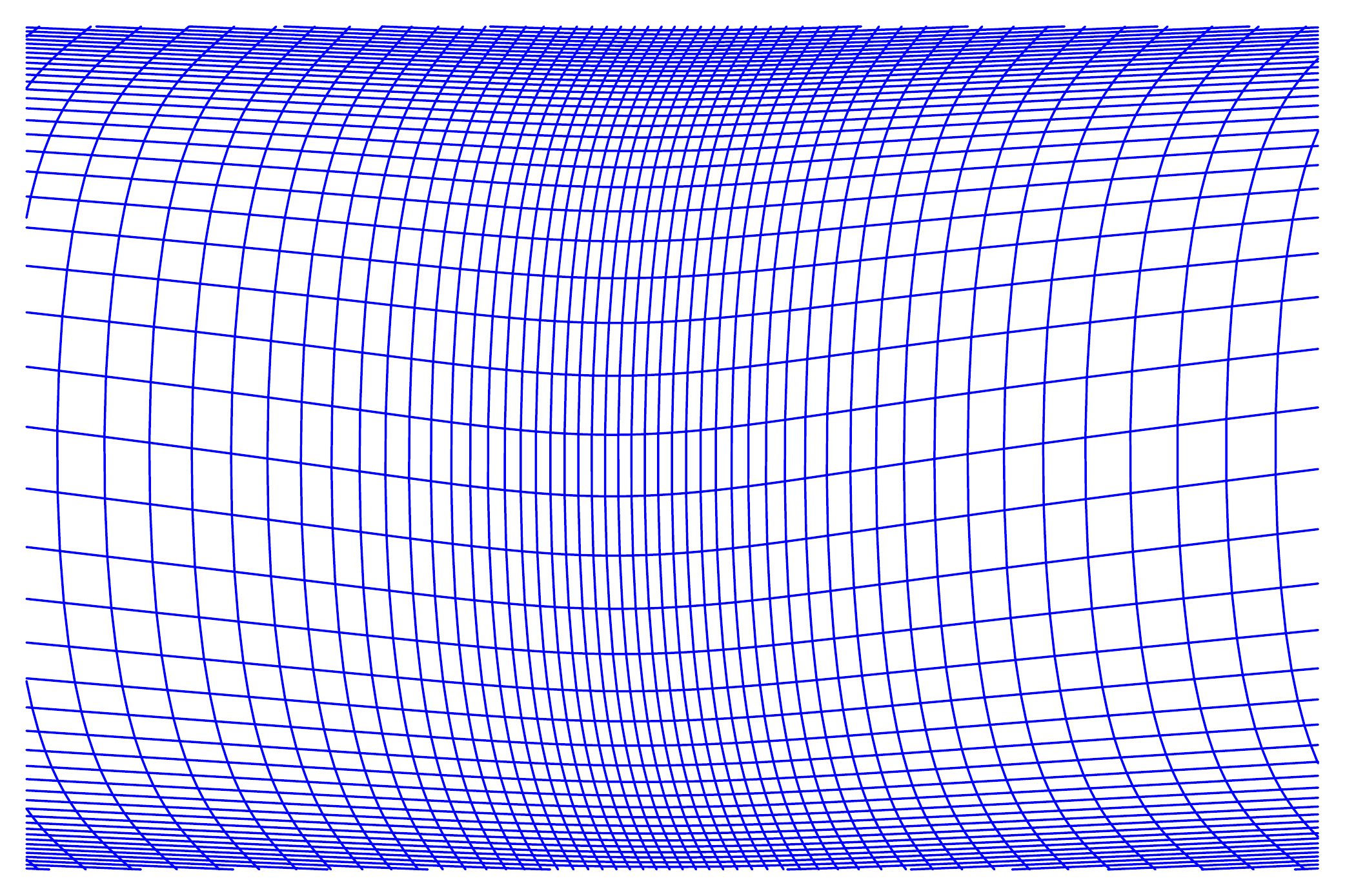}
         \caption{ Grid deformation illustrating the warp. }
         \label{fig:4d_h01_warp.pdf}
         \end{subfigure}
         \caption{The result of applying the shooting method with four  vector fields that do not form an algebra. The resulting deformation is a drastic improvement to the previous examples.} 
         \label{fig:4d}
\end{figure}

The dynamics of $u$ are in this case given by 
\begin{align*}
    & \dot u_1 =  5 u_3^2/32, \quad ~~\, \dot u_2 = -u_2 u_4^2, \\
    & \dot u_3 =  -u_1u_3/4, \quad  \dot u_4 =  5 u_2^2 /32.
\end{align*}
With all parameters as before and using the L-BFGS method, we obtain the moved landmarks and warp shown in \autoref{fig:4d}. 

We remark that while the computed deformation is by no means perfect, we obtain a much better deformation, and this was achieved by selecting vector fields whose components functions were trigonometric functions. 

\begin{figure}[!htb]
     \centering
     \begin{subfigure}[t]{0.48\textwidth}
         \centering
         \includegraphics[width=\textwidth]{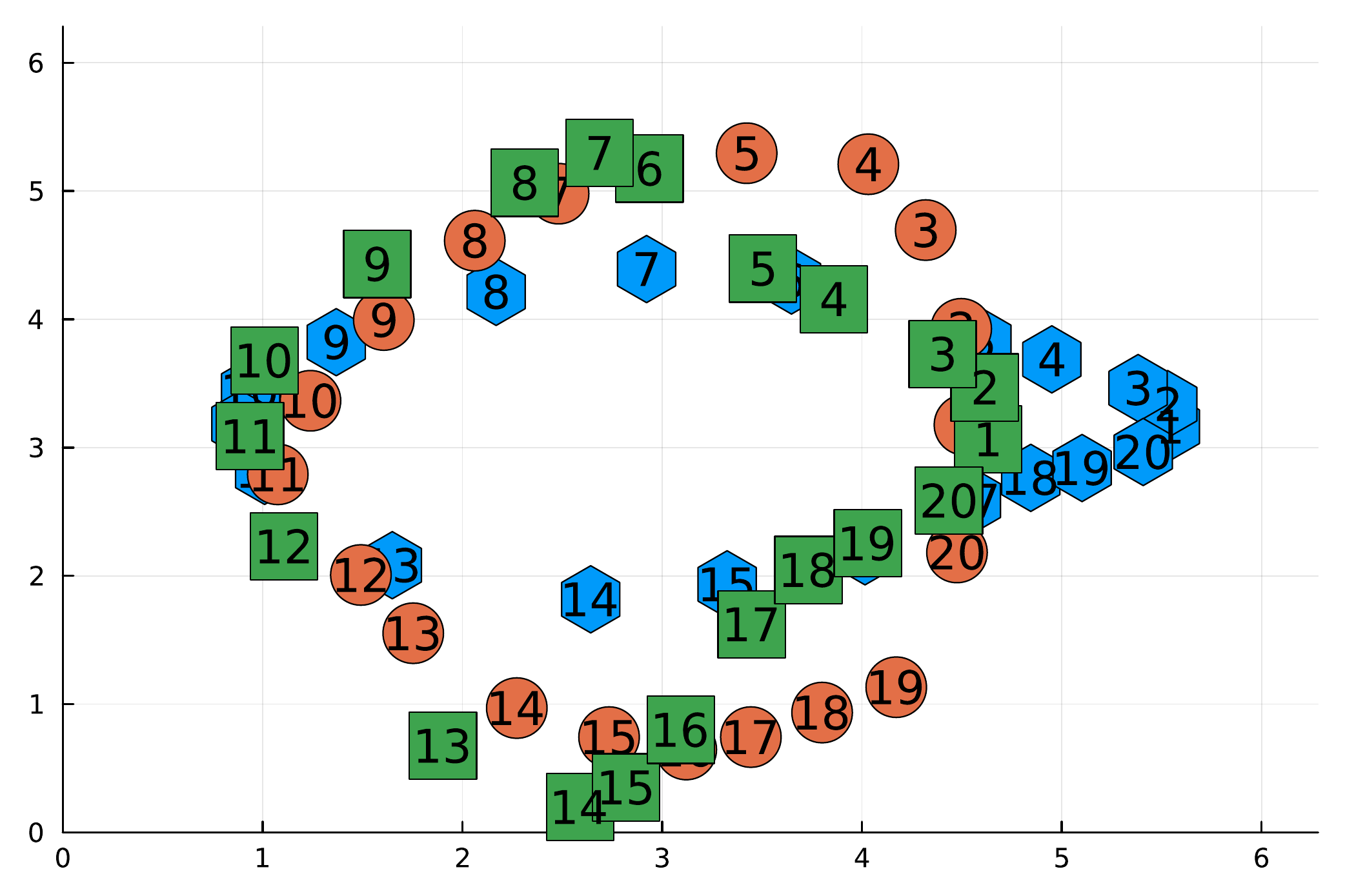}
         \caption{$h= 1$.}
         \label{fig:4d_result_h1}
         \end{subfigure}
         ~
         \begin{subfigure}[t]{0.48\textwidth}
         \centering
           \includegraphics[width=\textwidth]{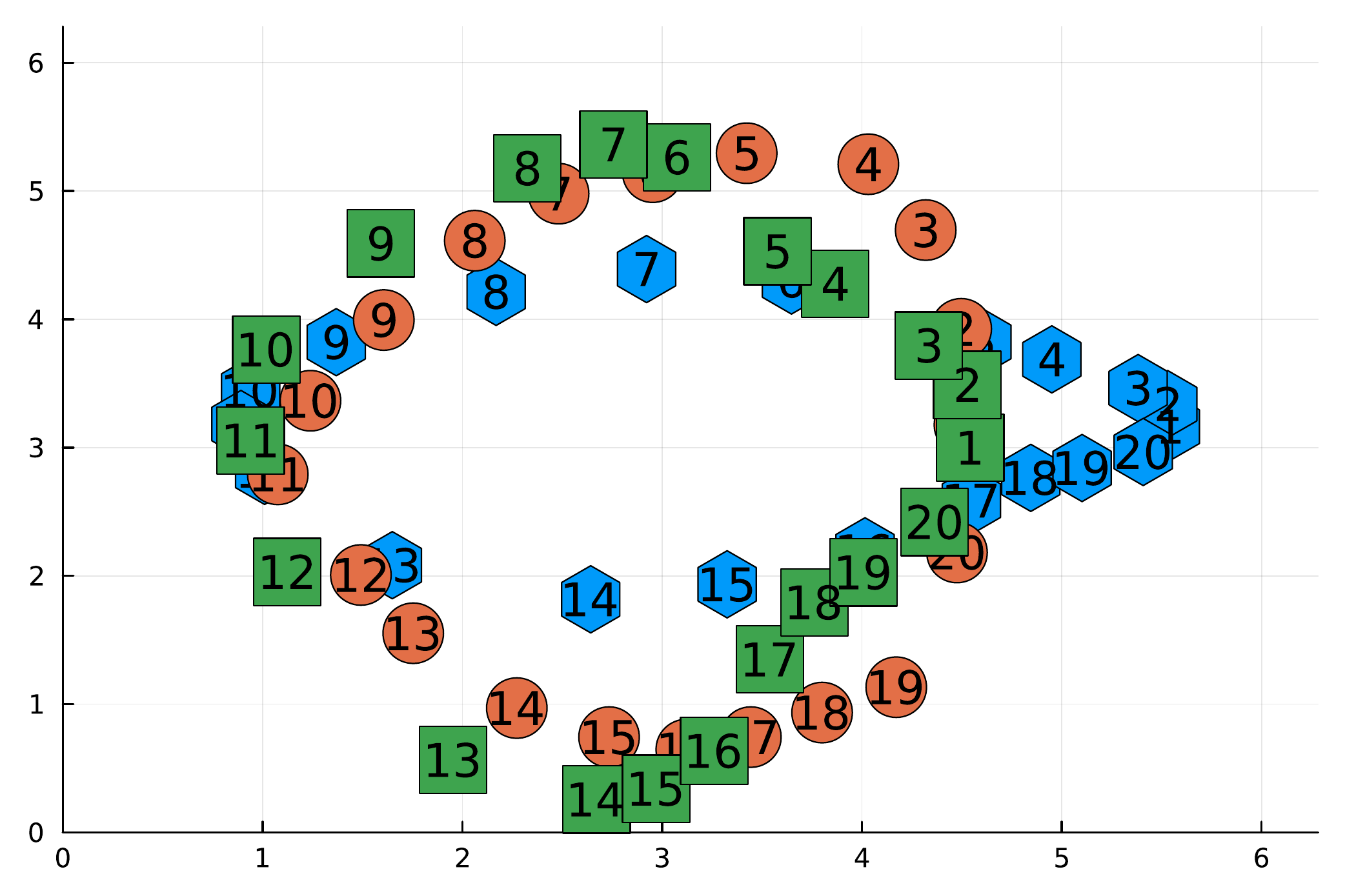}
         \caption{ $h=0.5$. }
         \label{fig:4d_result_h05}
         \end{subfigure}
         \caption{The effect of changing the integration step size $h$. The blue hexagons are the initial landmarks, the orange circles are the targets and the green squares are the deformed landmarks.
         Note that the matching improves as $h$ decreases. This is because the vector fields do not form an algebra.} 
         \label{fig:4d_hdiff}
\end{figure}

Contrary to the previous example, increasing or decreasing the number of time-steps does have a direct effect on the quality of the warp. 
Indeed, repeating the experiment with $h = 1$ and $h= 0.5$ yields the matching depicted in \autoref{fig:4d_hdiff}.
This is because the vector fields do not form an algebra, so taking more and more steps yields ever more iterated brackets. 
We observe that while the difference between $h=1$ and $h=0.1$ is quite drastic, the quality of the warp generally increases little for smaller $h$, at the expense of computational efficiency. 

\subsubsection{Vector fields given by truncated Fourier expansions}

In the previous experiment, we saw how trigonometric component functions yielded interesting dynamics. 
We continue this line of thought and select the vector fields to consist of eigenfunctions of the toroidal Laplace--Beltrami operator.

The control system is a linear combination of eigenfunctions, where the time-dependent weights are the control parameters. 
The problem of selecting the weights can therefore be seen as finding the best time-dependent vector field in the span of the (truncated) Fourier basis of eigenfunctions.  
While the number of basis functions included in each direction directly  influences how complicated vector fields we can approximate, we remark that already for a few terms, the iterated Lie brackets generate a large family of vector fields, and therefore we should be able to obtain quite complicated warps.

To illustrate, we select $1,\cos(x), \sin(x)$, $\cos(y), \sin(y)$ and in both the $x$ and $y$-direction, resulting in a total number of $10$ vector fields, so that $\cU = \R^{10}$.
The vector fields $X^1,\dots,X^{10}$ are explicitly
\begin{align*}
    &X^i = \psi_i(x,y) e_1 \text{  for } i=1,...,5\\
    &X^i = \psi_{i-5}(x,y) e_2  \text{  for } i=6,...,10
\end{align*}
where $\psi_i(x,y) = 1_{i=1,2} \cos((i-1)x)+ 1_{i=3} \sin(x)+ 1_{i=4} \cos(y)+ 1_{i=5} \sin(y)$.
Note that these vector fields are orthogonal, so $A(u)$ is invertible. 
The ODE governing the dynamics of $u$ is quadratic, as in the previous examples.
Therefore, we have local existence and uniqueness.

\begin{figure}[!htb]
     \centering
     \begin{subfigure}[t]{0.48\textwidth}
         \centering
         \includegraphics[width=\textwidth]{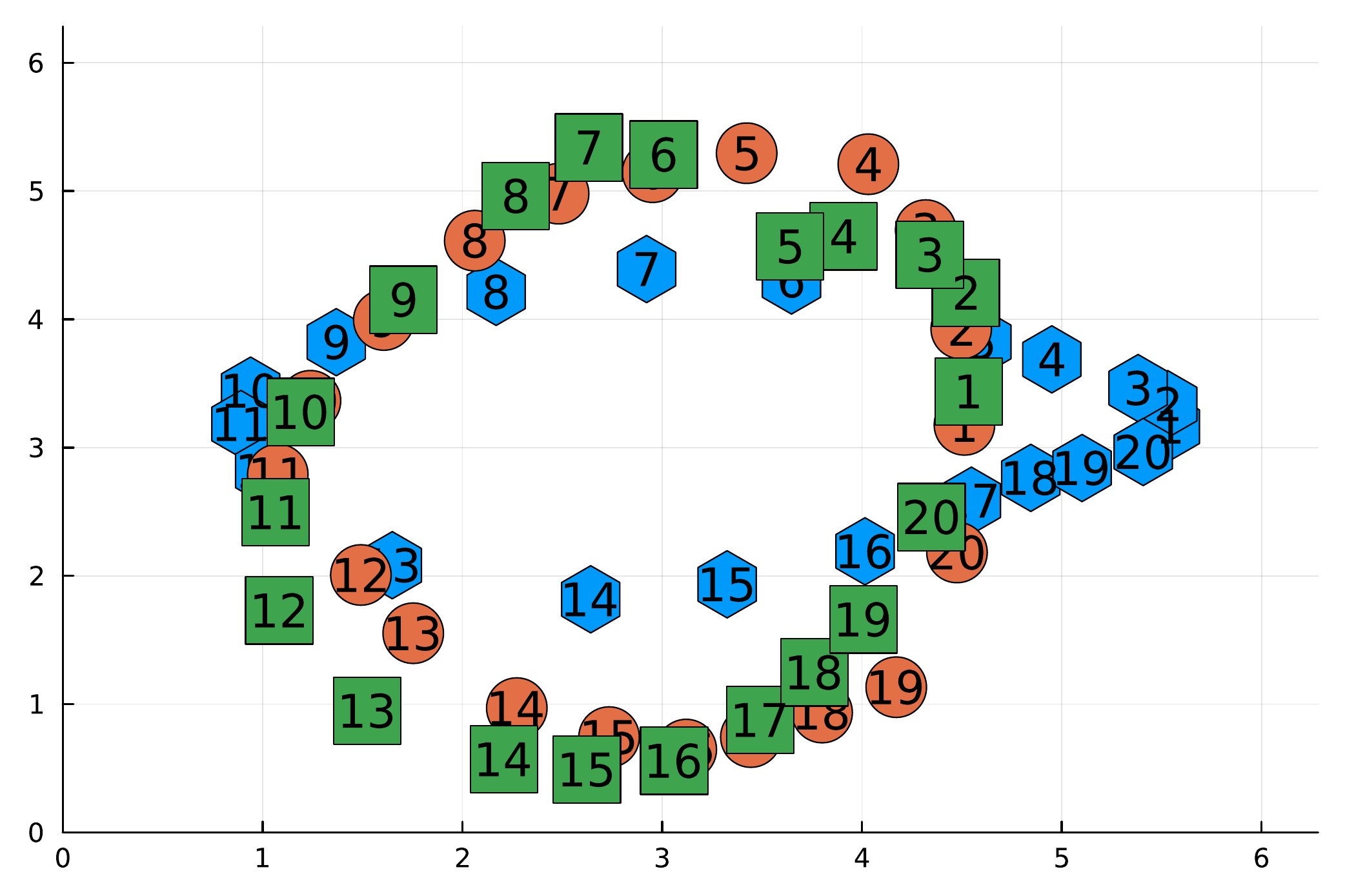}
         \caption{The blue hexagons are the initial landmarks, the orange circles are the targets and the green squares are the deformed landmarks. }
         \label{fig:result}
         \end{subfigure}
         ~
        \begin{subfigure}[t]{0.48\textwidth}
         \centering
          \includegraphics[width=\textwidth]{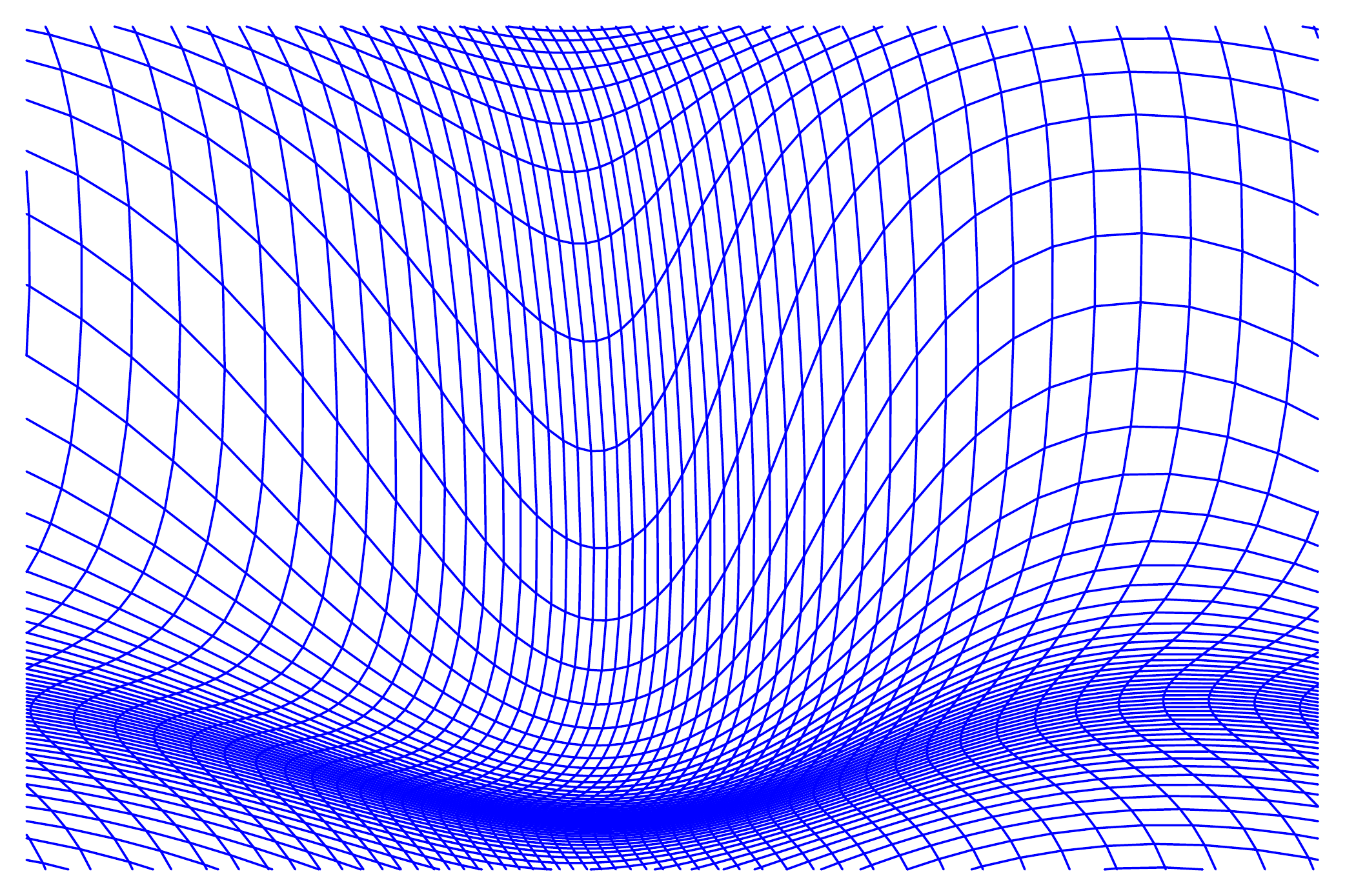}
         \caption{ Grid deformation illustrating the warp. 
         }
         \label{fig:warp}
         \end{subfigure}
         \caption{The result of applying the shooting method with more complicated vector fields, taken to build truncated Fourier expansions. }
         \label{fishres}
\end{figure}

We choose a regularization strength of $\sigma = 1/(400\pi^2)$.
With all other parameters and methods as before, the shooting method results in the moved landmarks in \autoref{fig:result}. 
In \autoref{fig:warp} the warp  is visualized as a grid deformation. 
Note that again, we do not obtain a perfect warp, but it is an improvement to only including $4$ vector fields. 
Further, including $10$ vector fields means that the family of obtained Lie brackets is rich, meaning that quite complicated warps should be available. 
The conclusion from the previous example, that increasing $h$ generally leads to worse matching, still holds in this case.

\subsection{Forward model example}

\begin{figure}[!htb]
     \centering
    \includegraphics[width = 0.7\linewidth]{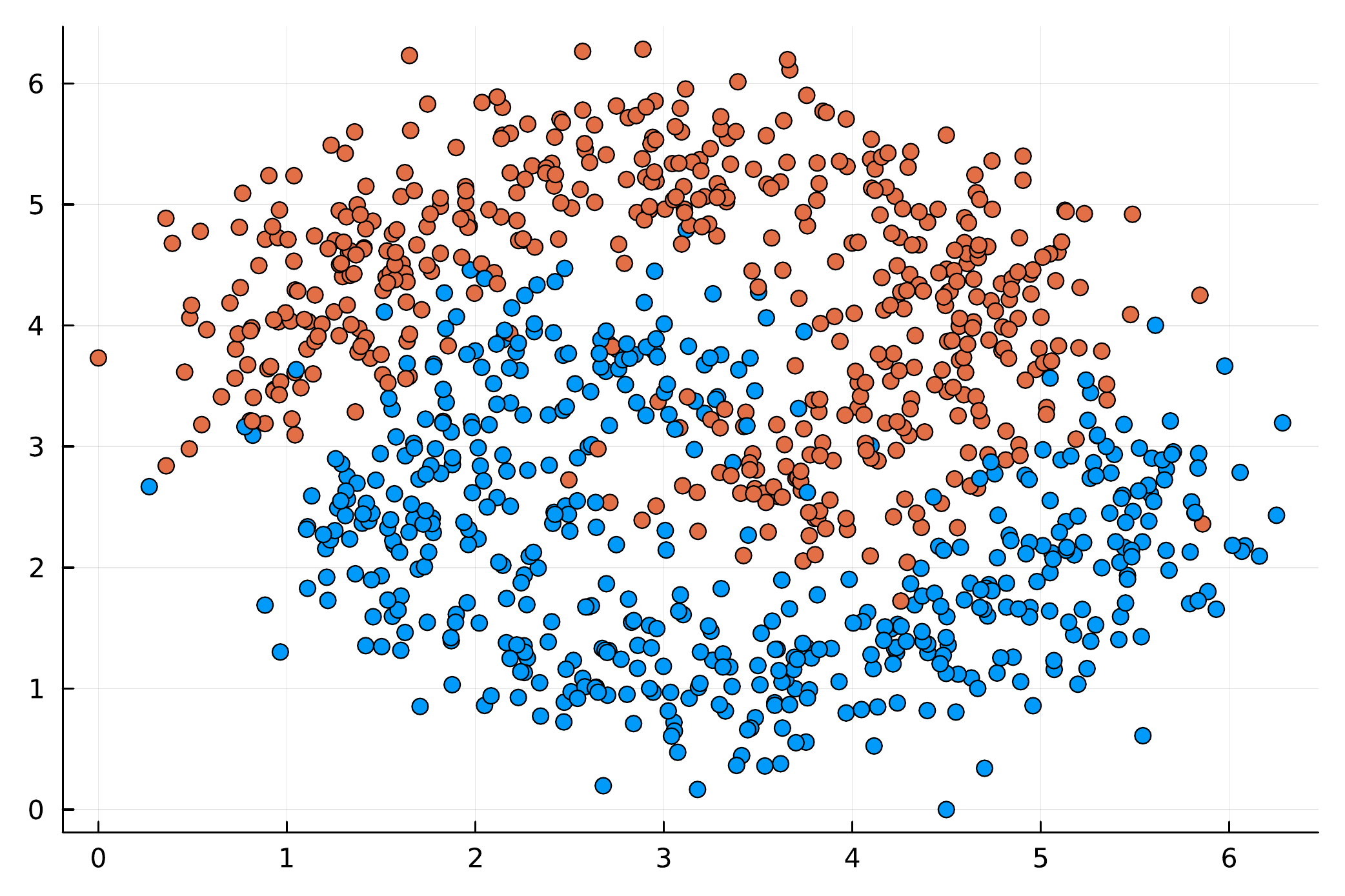}
         \caption{Two artificial datasets. The orange dots are classified as $1$, the blue dots as $0$.}
        \label{fig:spiral}
\end{figure}

To illustrate how the forward model can be used, we consider an artificial data set of initial landmarks consisting of two spiral arms.
We create targets by assigning $0$ to all landmarks in the lower spiral, and $1$ to all landmarks in the upper spiral. 
The data set is illustrated in \autoref{fig:spiral}. 
In this case, $N = [0,1]$ equipped with the usual metric. 
We select a forward model given by 
\begin{align*}
    h(x,y) = 16^{-(y-\pi)^2/\pi^2}. 
\end{align*}
This function is chosen so that points $(x,y) \in \R^2 /(2\pi \Z)^2$ where $y \in [\pi/2,3\pi/2]$ are mapped into $[1/2,1]$ and are thus closer to $1$ than to $0$, whereas points where $y \in [0,\pi/2] \cup [3\pi/2,2\pi]$ are mapped into $[0,1/2]$ and are closer to $0$. 
While we reiterate that this experiment is performed to illustrate the inclusion of a forward model, we remark that the setup is similar to a classification problem. 
Indeed, points in the band $\R^2 /(2\pi \Z)^2 \ni (x,y) \in [0,2\pi] \times (\pi/2,3\pi/2]$ are classified as a $1$, whereas the points outside the band are classified as $0$. 
The classifier $h(x)$ is chosen so that around $50\%$ of the points in a spiral are classified correctly.

The goal of the landmark deformation is therefore to unwrap the spiral, so that the lower spiral arm is moved outside the band whereas the upper spiral arm is moved inside. 
We run the experiment twice, varying the regularization strength to illustrate its impact on the diffeomorphic qualities of the found warp. 

\begin{figure}[!htb]
     \centering
     \begin{subfigure}[t]{0.48\textwidth}
         \centering
         \includegraphics[width=\textwidth]{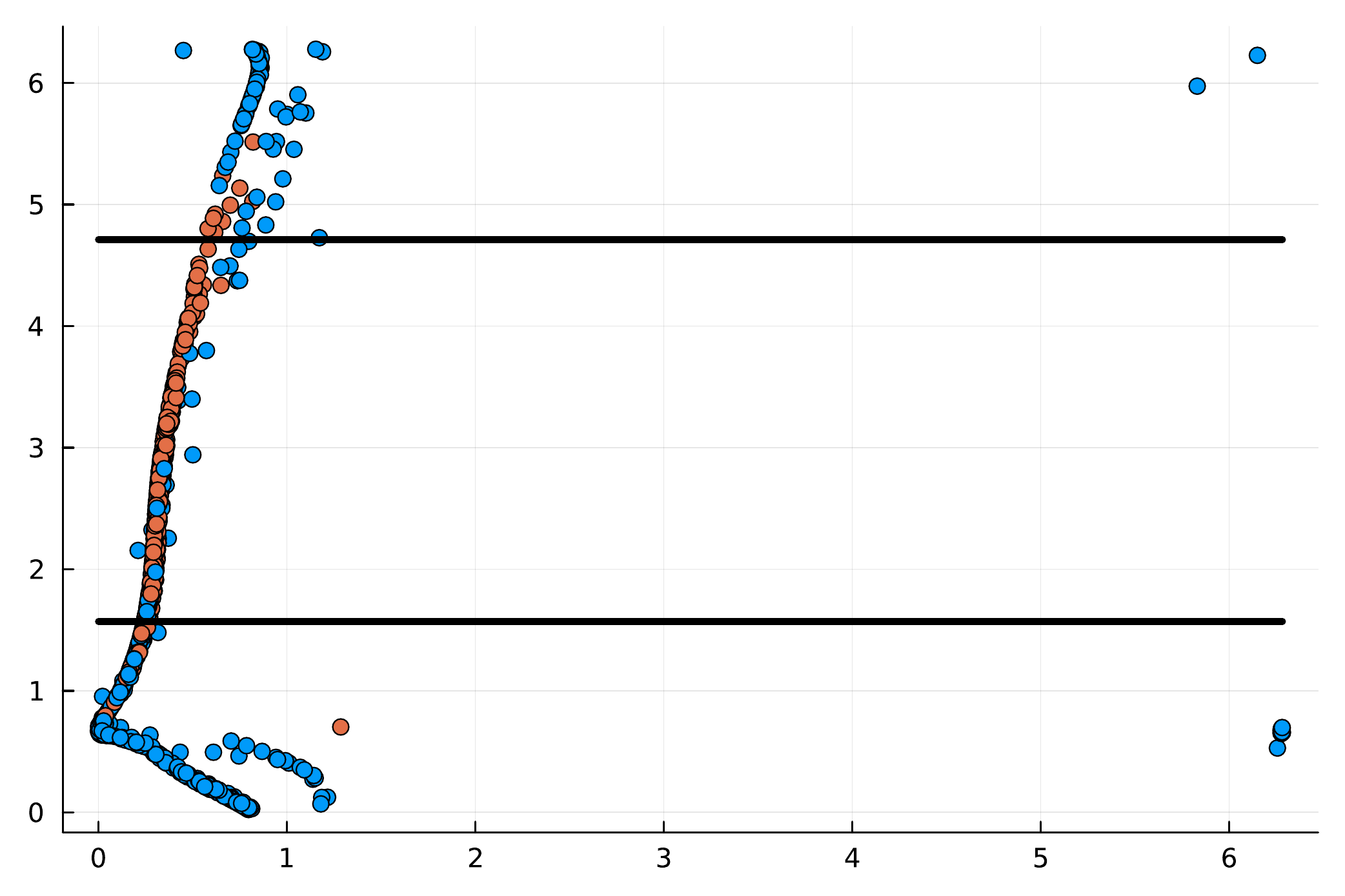}
         \caption{ The blue points are, in general, moved towards the line $y=\pi$, thus decreasing the distance to $1$ after applying the forward model. }
         \label{fig:result_lowreg}
         \end{subfigure}
         ~
         \begin{subfigure}[t]{0.48\textwidth}
         \centering
           \includegraphics[width=\textwidth]{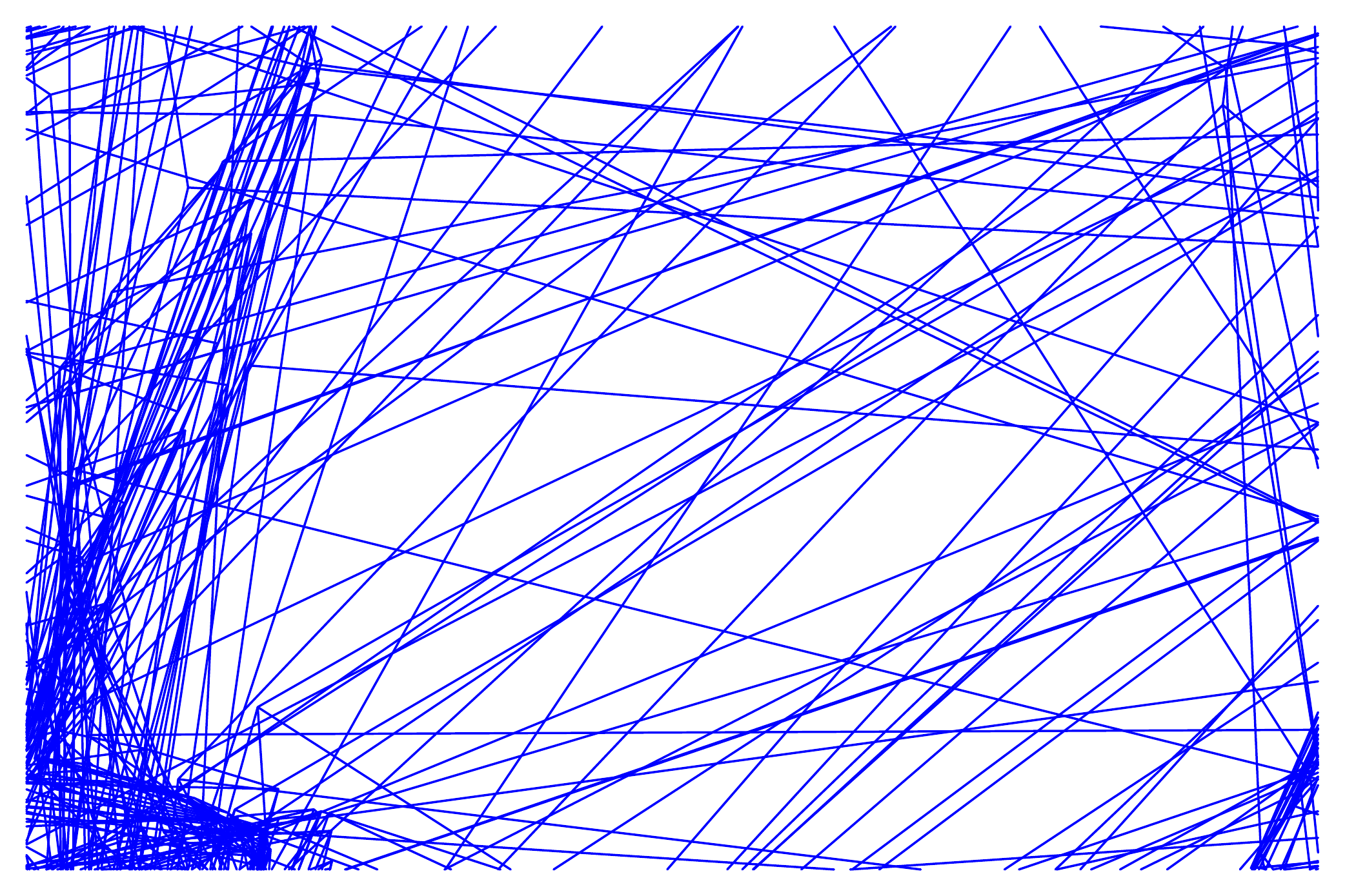}
         \caption{ Grid deformation illustrating the warp. Note that the warp is very ill-behaved.}
         \label{fig:warp_lowreg}
         \end{subfigure}
         \caption{The result of applying the shooting method to the spiral data set using a lower regularization parameter. 
         The classification boundary consists of the black horizontal lines: points between them are classified as $1$, and points outside them are classified as $0$.}
         \label{fig:res_spiral_lowreg}
\end{figure}

We first use a regularization strength of $\sigma = 1/(20\pi^2)$. 
The setting and parameter values are otherwise unchanged from the previous landmark matching experiment where $\cU = \R^{10}$ and the vector fields are built from the Laplace--Beltrami eigenfunctions. 
The shooting method with the L-BFGS algorithm to optimize $u_0$, initialized as the zero vector, yields the results in \autoref{fig:res_spiral_lowreg}. 
The warp in \autoref{fig:warp_lowreg} is clearly not well-behaved from a diffeomorphic perspective. 
However, after applying the warp, almost $87\%$ of the points are classified correctly. 
We now increase the regularization to see if the diffeomorphic properties of the warp can be improved.  

\begin{figure}[!htb]
     \centering
     \begin{subfigure}[t]{0.48\textwidth}
         \centering
         \includegraphics[width=\textwidth]{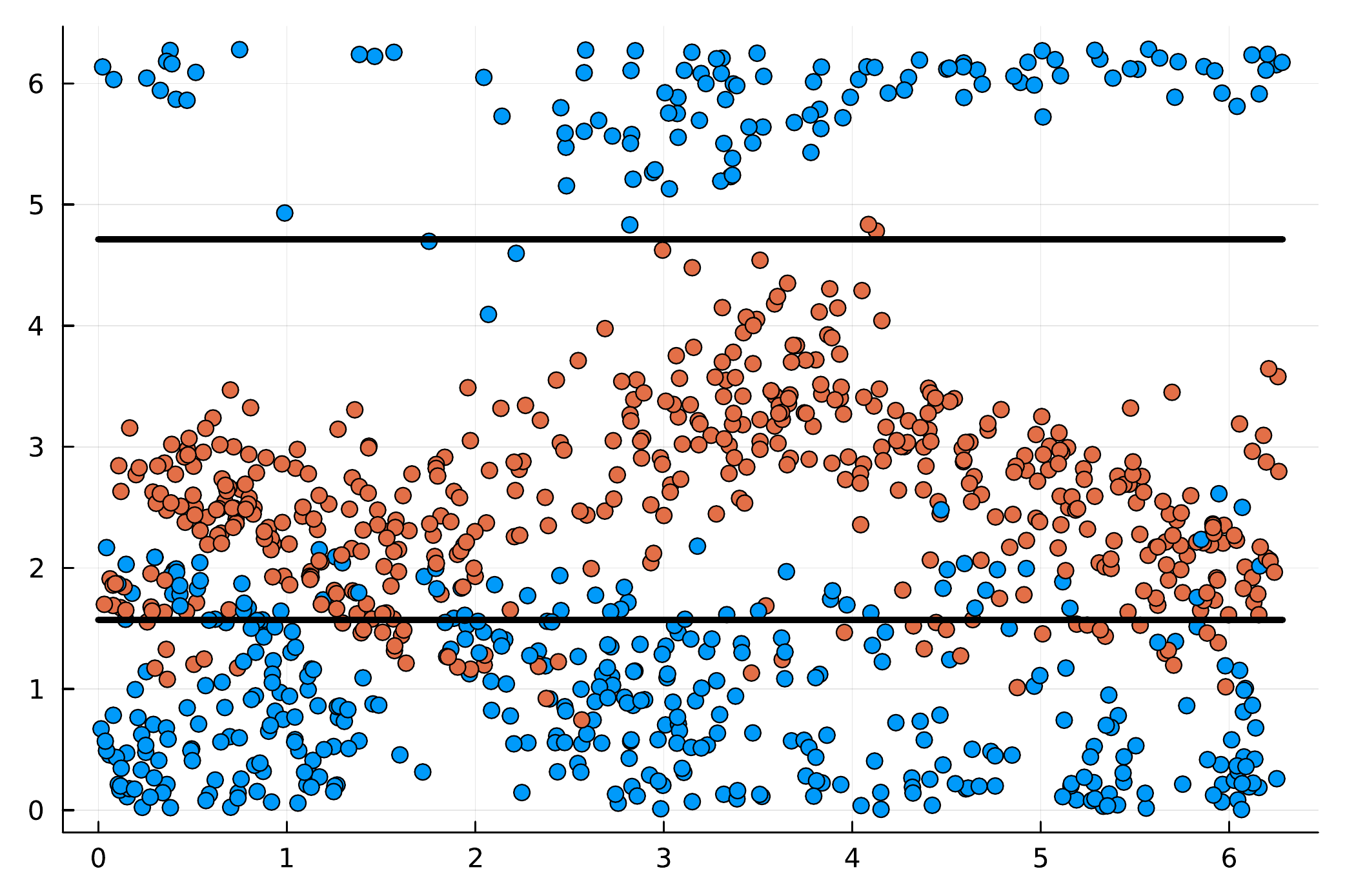}
         \caption{ The blue points are, in general, moved towards the line $y=\pi$, thus decreasing the distance to $1$ after applying the forward model. }
         \label{fig:result_highreg}
         \end{subfigure}
         ~
         \begin{subfigure}[t]{0.48\textwidth}
         \centering
           \includegraphics[width=\textwidth]{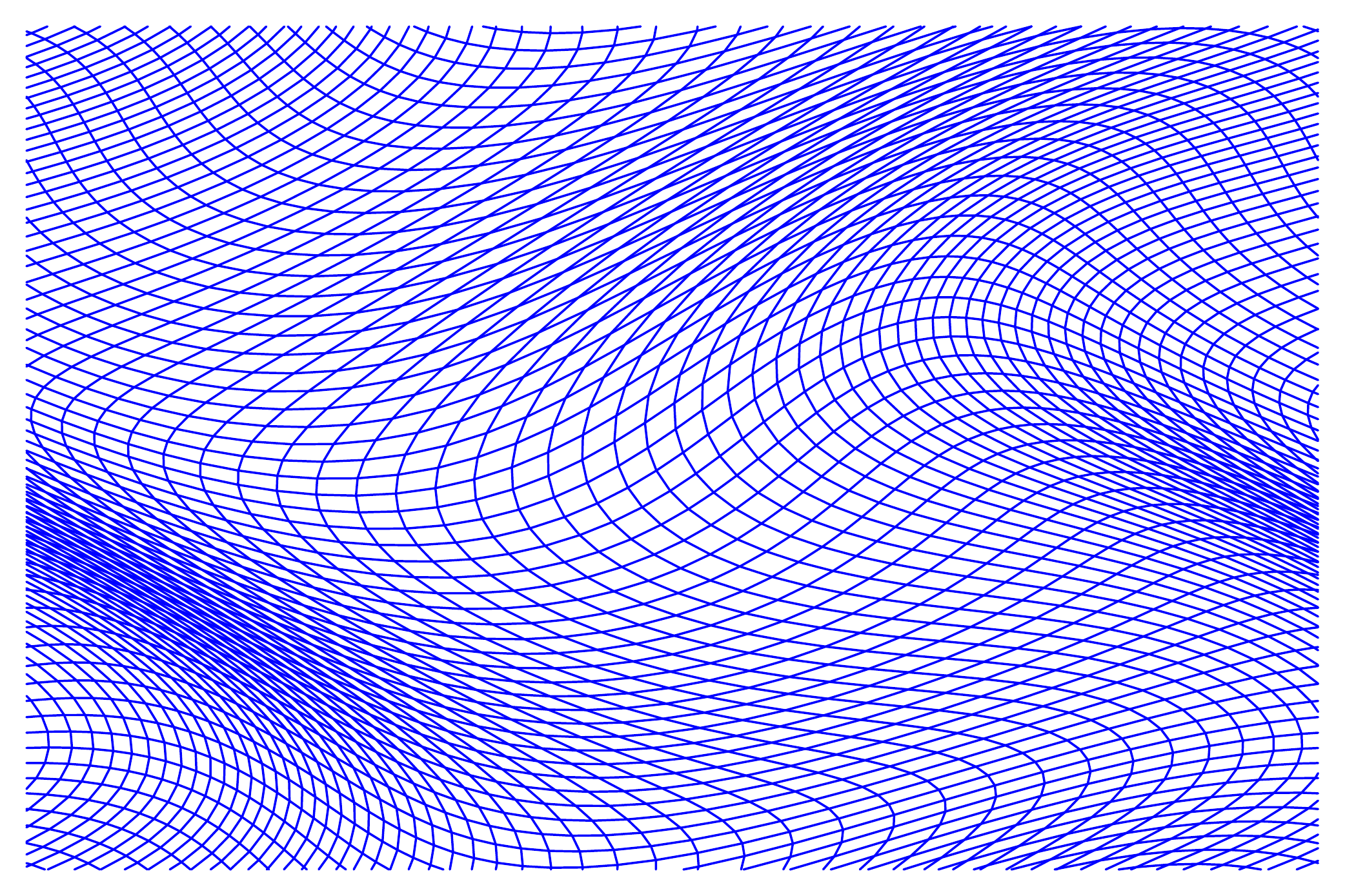}
         \caption{ Grid deformation illustrating the warp.}
         \label{fig:warp_highreg}
         \end{subfigure}
         \caption{The result of applying the shooting method to the spiral data set using a higher regularization parameter. 
         The classification boundary consists of the black horizontal lines: points between them are classified as $1$, and points outside them are classified as $0$.}
         \label{fig:res_spiral_highreg}
\end{figure}

We set the regularization strength to $1/(4\pi^2)$ and obtain the results in \autoref{fig:res_spiral_highreg}. 
In this case, the found warp is better behaved. 
Further, the classification performance is comparable to the low-regularized case: about $87\%$ of the points are classified correctly after the warp is applied.

This warp, shown in \autoref{fig:warp_highreg} can be compared with the warp depicted in \autoref{fig:warp_lowreg}, generated with a small regularization parameter, which is ill-behaved.
Indeed, the result in \autoref{fig:warp_lowreg} is not \emph{computationally stable}~\citep{Webb}: nearby points do not necessarily remain close after the warp is applied.
The computational stability is much better with the increased regularization parameter, as depicted in \autoref{fig:warp_highreg}.  
We conclude that one must balance the regularization strength to find warps that are useful but not too erratic.

We observed in the two examples that the classification performance was comparable for the two considered regularization strengths.
To further investigate this, we run the same experiment for various values of the regularization strength $\sigma$ and compare the classification performance on the points in the original spiral as well as on $10\,000$ new points generated in the same way as the original points.
 We remark that the regularization strength cannot be taken too low, as this results in the optimization becoming unstable. 
This is related to the breakdown of the diffeomorphic properties as the regularization decreases.

\begin{table}[h]
\centering
\begin{tabular}{@{}c c c@{}}
\toprule
\textbf{$\sigma$ (times $1/(4\pi^2)$)} & \textbf{Correct $\%$} & \textbf{Correct $\%$, new data} \\
\midrule 
$0.2$              & $86.5$     & $83.2$ \\ 
$1.0$   & $87.8$            & $86.9$ \\ 
$10.0$     & $83.7$  & $82.8$ \\ 
$100.0$   & $78.8$        & $75.4$ \\ 
$1000.0$   & $60.2$        & $53.5$ \\ 
\bottomrule
\end{tabular}
\caption{The impact of regularization strength on classification performance. }\label{tbl:warp_quality}
\end{table}

\autoref{tbl:warp_quality} suggests that there is a balance to strike between the diffeomorphic qualities of the warp, i.e., the regularization strength, and the classification performance. 
 Indeed, when the diffeomorphic qualities increase slightly, from the lowest considered regularization strength, the classification performance increases slightly.
However, increasing the regularization too much results in poor classification performance.

Further, comparing the performance on the original data to the new data,
we see that for a balanced choice of regularization strength (i.e., $\sigma = 1/(4\pi^2)$ or $\sigma = 10/(4\pi^2)$), the 
 performance on new data is not significantly decreased. 
 This is not the case for lower and higher regularization strengths.
 
This outcome suggests a potentially interesting future study on the relationship between diffeomorphic qualities and classification performance and generalization.

\section{The ResNet connection}
\label{sec:resnet}
In this section, we discuss the connection between ResNets and sub-Riemannian landmark matching. 
We begin with  a brief overview of ResNets and their interpretation as temporal discretizations of continuous time optimal control problems.
It is this interpretation that allows for a comparison with sub-Riemannian landmark matching.

\begin{figure}[!htb]
     \centering
     \begin{subfigure}[t]{0.45\linewidth}
        \centering
        \scalebox{0.7}{
        \tikzset{every picture/.style={line width=0.75pt}} 
        
        \begin{tikzpicture}[x=0.75pt,y=0.75pt,yscale=-1,xscale=1]
      
        \draw  [fill={rgb, 255:red, 239; green, 231; blue, 231 }  ,fill opacity=1 ] (275,120) .. controls (275,112.27) and (281.27,106) .. (289,106) -- (331,106) .. controls (338.73,106) and (345,112.27) .. (345,120) -- (345,190) .. controls (345,197.73) and (338.73,204) .. (331,204) -- (289,204) .. controls (281.27,204) and (275,197.73) .. (275,190) -- cycle ;
      
         \draw [line width=1.5]    (230,155) -- (271,155) ;
        \draw [shift={(275,155)}, rotate = 180] [fill={rgb, 255:red, 0; green, 0; blue, 0 }  ][line width=0.08]  [draw opacity=0] (8.13,-3.9) -- (0,0) -- (8.13,3.9) -- cycle    ;
     
        \draw [line width=1.5]    (346,155) -- (400,155) ;
        \draw [shift={(404,155)}, rotate = 180] [fill={rgb, 255:red, 0; green, 0; blue, 0 }  ][line width=0.08]  [draw opacity=0] (8.13,-3.9) -- (0,0) -- (8.13,3.9) -- cycle    ;
        \draw [shift={(346,155)}, rotate = 180] [color={rgb, 255:red, 0; green, 0; blue, 0 }  ][line width=1.5]    (0,4.7) -- (0,-4.7)   ;
        
        \draw (240,118.4) node [anchor=north west][inner sep=0.75pt]    {$Y^{[i]}$};
        \draw (376,116.4) node [anchor=north west][inner sep=0.75pt]    {$Y^{[i+1]}$};
        \end{tikzpicture}
        }
        \caption{A layer in a feed forward neural network.}
                \label{fig:dnn}    \end{subfigure}
         ~
    \begin{subfigure}[t]{0.45\linewidth}
        \centering
        \scalebox{0.7}{
        \tikzset{every picture/.style={line width=0.75pt}} 
        
        \begin{tikzpicture}[x=0.75pt,y=0.75pt,yscale=-1,xscale=1]
            \draw  [fill={rgb, 255:red, 239; green, 231; blue, 231 }  ,fill opacity=1 ] (275,120) .. controls (275,112.27) and (281.27,106) .. (289,106) -- (331,106) .. controls (338.73,106) and (345,112.27) .. (345,120) -- (345,190) .. controls (345,197.73) and (338.73,204) .. (331,204) -- (289,204) .. controls (281.27,204) and (275,197.73) .. (275,190) -- cycle ;
       
            \draw [line width=1.5]    (183,155) -- (271,155) ;
            \draw [shift={(275,155)}, rotate = 180] [fill={rgb, 255:red, 0; green, 0; blue, 0 }  ][line width=0.08]  [draw opacity=0] (8.13,-3.9) -- (0,0) -- (8.13,3.9) -- cycle    ;
            \draw [line width=1.5]    (474,155) -- (502,155) ;
            \draw [shift={(506,155)}, rotate = 180] [fill={rgb, 255:red, 0; green, 0; blue, 0 }  ][line width=0.08]  [draw opacity=0] (8.13,-3.9) -- (0,0) -- (8.13,3.9) -- cycle    ;
            \draw [shift={(474,155)}, rotate = 180] [color={rgb, 255:red, 0; green, 0; blue, 0 }  ][line width=1.5]    (0,4.7) -- (0,-4.7)   ;
            \draw [line width=1.5]    (345,156) -- (433,156) ;
            \draw [shift={(437,156)}, rotate = 180] [fill={rgb, 255:red, 0; green, 0; blue, 0 }  ][line width=0.08]  [draw opacity=0] (8.13,-3.9) -- (0,0) -- (8.13,3.9) -- cycle    ;
            \draw [shift={(345,156)}, rotate = 180] [color={rgb, 255:red, 0; green, 0; blue, 0 }  ][line width=1.5]    (0,4.7) -- (0,-4.7)   ;
            \draw   (437,155) .. controls (437,144.92) and (445.28,136.75) .. (455.5,136.75) .. controls (465.72,136.75) and (474,144.92) .. (474,155) .. controls (474,165.08) and (465.72,173.25) .. (455.5,173.25) .. controls (445.28,173.25) and (437,165.08) .. (437,155) -- cycle ; \draw   (442.42,142.1) -- (468.58,167.9) ; \draw   (468.58,142.1) -- (442.42,167.9) ;
            \draw    (229,155) .. controls (229,-5.69) and (455.72,16.27) .. (455.52,134.95) ;
            \draw [shift={(455.5,136.75)}, rotate = 271.19] [fill={rgb, 255:red, 0; green, 0; blue, 0 }  ][line width=0.08]  [draw opacity=0] (8.93,-4.29) -- (0,0) -- (8.93,4.29) -- cycle    ;
            
            \draw (191,118.4) node [anchor=north west][inner sep=0.75pt]    {$Y^{[ j]}$};
            \draw (460,110.4) node [anchor=north west][inner sep=0.75pt]    {$Y^{[ j+1]}$};
            \draw (349,112.4) node [anchor=north west][inner sep=0.75pt]  [font=\footnotesize]  {$f\left( K^{[ j]} Y^{[ j]} +b^{[ j]}\right)$};
            \draw (442,49.4) node [anchor=north west][inner sep=0.75pt]    {$Y^{[ j]}$};
            \end{tikzpicture}
            }
            \caption{ResNet layer. }
            \label{fig:resnet}   
        \end{subfigure}
        \caption{Schematic illustration of layers in two different neural network architectures.}
        \label{fig:nets}
\end{figure}
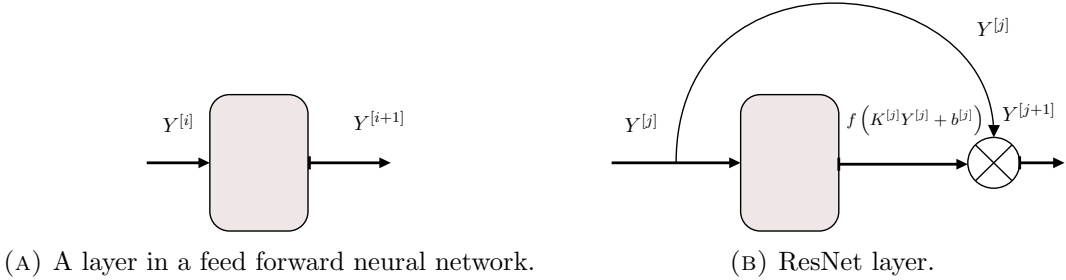

Given input data $X=(X_1,...,X_n)=Y^{[0]}=(Y_1^{[0]},...,Y_n^{[0]}) \in \R^{d}$ and output data $Z=(Z_1,...,Z_n)$, the \emph{supervised learning problem} is to find a set of transformation parameters $u$ and classification weights $C=(W,\mu)\in R^{d \times d}\times \R^{d}$ so that the function 
\begin{align}
\cL\left(h(Y^{[N]},W),Z\right)+\cR\left(u\right),
    \label{eq:resnetfunc}
\end{align}
is minimized for some loss function $\cL$, classifier $h$ and regularizer $\cR$. 
Here, $Y^{[N]}$ is given by  warping the initial data. 
If neural networks are used to parametrize the warp, the so-called \emph{neural architecture} is determined by the recursive relation between $Y^{[i]}$ and $Y^{[i+1]}$.
In the case of a standard feed-forward neural network, see \autoref{fig:dnn} for a depiction of a layer, this relation is typically given by 
\begin{align*}
    Y_{k}^{[i+1]}=f(Y_k^{[i]},u^{[i]}),
\end{align*} where $u^{[i]}$ are weights and biases of the $(i+1)$th layer and $f$ denotes some suitable non-linear function, often called the \emph{activation function}.
In the case of ResNets, a \emph{skip connection} is added between the input layer and its output. A layer is depicted in \autoref{fig:resnet}.
In this case, the relation between $Y^{[i+1]}$ and $Y^{[i]}$ will be given by  
\begin{align*}
    Y_{k}^{[i+1]}=Y_k^{[i]}+hf(Y_k^{[i]},u^{[i]}).
\end{align*}
Skip connections are used to alleviate training issues stemming from the \emph{vanishing gradient problem} \citep{He2015}.
	
Combining  \autoref{eq:resnetfunc} with the skip connection architecture,  we see that a supervised ResNet training problem can be written as an optimization problem,
\begin{align}
    &\!\min_{Y,u,C} \quad \cL\left(h\left(Y^{[N]},C\right),Z\right)+\cR\left(u\right) \label{eq:target},\\
    &\text{s.t.  } \quad Y_{k}^{[i+1]}=Y_k^{[i]}+hf(Y_k^{[i]},u^{[i]}),~i=0,1,2,...,N-1, Y^{[0]}_i=X_i,  k=1,...,n \label{eq:condition}.
\end{align}

It is possible to consider the neural architecture of \autoref{eq:condition} as the forward Euler discretization of the ODE $\dot y_i = f(y_i,u)$  \citep{Celledoni2020,Li2017}.
The ResNet training problem can therefore be considered as a discrete time version of the optimal control problem 
\begin{align}
    &\!\min_{y,u,C} \quad \cL\left(h\left(y(T),u\right),Z\right)+\cR\left(u\right) \label{eq:targetcont}\\
        &\text{s.t.  } \quad \dot y_i = f(y_i,u(t)), t\in [0,1], ~y_i(0)=X_i  \label{eq:conditioncont}.
\end{align}

Let us now compare with the problem in \autoref{eq:targetvlm_nonformal}-\ref{eq:constvlm_nonformal}. 
Considering the case where $\cL = \sum_{i=1}^n d_N^2(h(y_i(1)),c_i)$ it is clear that we can interpret sub-Riemannian landmark matching as the continuous-time optimization problem stemming from a deep learning problem with a ResNet architecture.

In the ResNet case, $u_t$ will be given by weight and bias curves $(K(t),b(t))$ where $$K:[0,1] \to \R^{d \times d}$$
 and  $b\colon[0,1] \to \R^{d}$. 
Letting $f(t,u)(\cdot)= f(K(t)\cdot+b(t))$, we have a map from the set of controls $\cU= \R^{d\times d} \times \R^d$ into a set of vector fields on $\R^d$. 
In a ResNet, the forward model will be the final projection or classification function projecting the landmarks, in this case being images or input data, down to the space of labels.  

Assuming that the classification parameters and number of layers are fixed, we note that we are in the setting of \autoref{constrained_landmarks} by considering the map given by  $f(\cdot,u)$,  as this operator will map from spaces of time-dependent controls to spaces of time-dependent vector fields under some technical continuity and measurability assumptions.

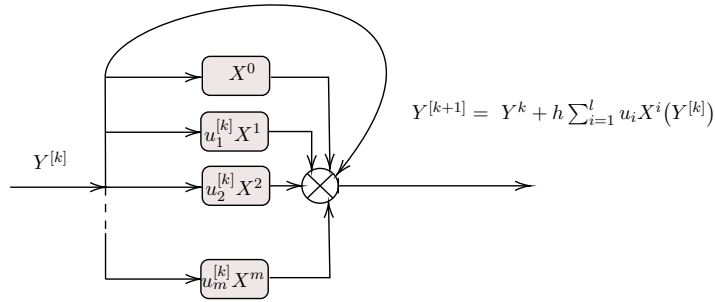
\begin{figure}[ht]
			\centering

\tikzset{every picture/.style={line width=0.75pt}} 
\scalebox{0.7}{
\begin{tikzpicture}[x=0.75pt,y=0.75pt,yscale=-1,xscale=1]

\draw    (108,150) -- (174,150) ;
\draw [shift={(176,150)}, rotate = 180] [color={rgb, 255:red, 0; green, 0; blue, 0 }  ][line width=0.75]    (10.93,-3.29) .. controls (6.95,-1.4) and (3.31,-0.3) .. (0,0) .. controls (3.31,0.3) and (6.95,1.4) .. (10.93,3.29)   ;
\draw    (175.75,71) -- (176,150) ;
\draw    (175.75,71) -- (241.75,71) ;
\draw [shift={(243.75,71)}, rotate = 180] [color={rgb, 255:red, 0; green, 0; blue, 0 }  ][line width=0.75]    (10.93,-3.29) .. controls (6.95,-1.4) and (3.31,-0.3) .. (0,0) .. controls (3.31,0.3) and (6.95,1.4) .. (10.93,3.29)   ;
\draw    (175.88,110.5) -- (200.88,110.5) -- (241.88,110.5) ;
\draw [shift={(243.88,110.5)}, rotate = 180] [color={rgb, 255:red, 0; green, 0; blue, 0 }  ][line width=0.75]    (10.93,-3.29) .. controls (6.95,-1.4) and (3.31,-0.3) .. (0,0) .. controls (3.31,0.3) and (6.95,1.4) .. (10.93,3.29)   ;
\draw    (176,150) -- (201,150) -- (242,150) ;
\draw [shift={(244,150)}, rotate = 180] [color={rgb, 255:red, 0; green, 0; blue, 0 }  ][line width=0.75]    (10.93,-3.29) .. controls (6.95,-1.4) and (3.31,-0.3) .. (0,0) .. controls (3.31,0.3) and (6.95,1.4) .. (10.93,3.29)   ;
\draw  [dash pattern={on 4.5pt off 4.5pt}]  (176,150) -- (176,183) ;
\draw    (176,216) -- (242,216) ;
\draw [shift={(244,216)}, rotate = 180] [color={rgb, 255:red, 0; green, 0; blue, 0 }  ][line width=0.75]    (10.93,-3.29) .. controls (6.95,-1.4) and (3.31,-0.3) .. (0,0) .. controls (3.31,0.3) and (6.95,1.4) .. (10.93,3.29)   ;
\draw    (176,183) -- (176,216) ;
\draw  [fill={rgb, 255:red, 239; green, 231; blue, 231 }  ,fill opacity=1 ] (245,61.6) .. controls (245,58.51) and (247.51,56) .. (250.6,56) -- (287.4,56) .. controls (290.49,56) and (293,58.51) .. (293,61.6) -- (293,78.4) .. controls (293,81.49) and (290.49,84) .. (287.4,84) -- (250.6,84) .. controls (247.51,84) and (245,81.49) .. (245,78.4) -- cycle ;
\draw    (175.75,71) .. controls (176,8) and (477,-26.5) .. (341,141) ;
\draw [shift={(341,141)}, rotate = 309.07] [color={rgb, 255:red, 0; green, 0; blue, 0 }  ][line width=0.75]    (10.93,-3.29) .. controls (6.95,-1.4) and (3.31,-0.3) .. (0,0) .. controls (3.31,0.3) and (6.95,1.4) .. (10.93,3.29)   ;
\draw   (316,149) .. controls (316,142.1) and (321.82,136.5) .. (329,136.5) .. controls (336.18,136.5) and (342,142.1) .. (342,149) .. controls (342,155.9) and (336.18,161.5) .. (329,161.5) .. controls (321.82,161.5) and (316,155.9) .. (316,149) -- cycle ; \draw   (319.81,140.16) -- (338.19,157.84) ; \draw   (338.19,140.16) -- (319.81,157.84) ;
\draw    (294,71) -- (335,71) -- (335.97,135) ;
\draw [shift={(336,137)}, rotate = 269.13] [color={rgb, 255:red, 0; green, 0; blue, 0 }  ][line width=0.75]    (10.93,-3.29) .. controls (6.95,-1.4) and (3.31,-0.3) .. (0,0) .. controls (3.31,0.3) and (6.95,1.4) .. (10.93,3.29)   ;
\draw    (292,110) -- (323,110) -- (323,136) ;
\draw [shift={(323,138)}, rotate = 270] [color={rgb, 255:red, 0; green, 0; blue, 0 }  ][line width=0.75]    (10.93,-3.29) .. controls (6.95,-1.4) and (3.31,-0.3) .. (0,0) .. controls (3.31,0.3) and (6.95,1.4) .. (10.93,3.29)   ;
\draw    (293,215) -- (336,215) -- (335.04,163) ;
\draw [shift={(335,161)}, rotate = 448.94] [color={rgb, 255:red, 0; green, 0; blue, 0 }  ][line width=0.75]    (10.93,-3.29) .. controls (6.95,-1.4) and (3.31,-0.3) .. (0,0) .. controls (3.31,0.3) and (6.95,1.4) .. (10.93,3.29)   ;
\draw    (342,149) -- (477,149) ;
\draw [shift={(479,149)}, rotate = 180] [color={rgb, 255:red, 0; green, 0; blue, 0 }  ][line width=0.75]    (10.93,-3.29) .. controls (6.95,-1.4) and (3.31,-0.3) .. (0,0) .. controls (3.31,0.3) and (6.95,1.4) .. (10.93,3.29)   ;
\draw [shift={(342,149)}, rotate = 180] [color={rgb, 255:red, 0; green, 0; blue, 0 }  ][line width=0.75]    (0,5.59) -- (0,-5.59)   ;
\draw  [fill={rgb, 255:red, 239; green, 231; blue, 231 }  ,fill opacity=1 ] (244,101.6) .. controls (244,98.51) and (246.51,96) .. (249.6,96) -- (286.4,96) .. controls (289.49,96) and (292,98.51) .. (292,101.6) -- (292,118.4) .. controls (292,121.49) and (289.49,124) .. (286.4,124) -- (249.6,124) .. controls (246.51,124) and (244,121.49) .. (244,118.4) -- cycle ;
\draw  [fill={rgb, 255:red, 239; green, 231; blue, 231 }  ,fill opacity=1 ] (245,140.6) .. controls (245,137.51) and (247.51,135) .. (250.6,135) -- (287.4,135) .. controls (290.49,135) and (293,137.51) .. (293,140.6) -- (293,157.4) .. controls (293,160.49) and (290.49,163) .. (287.4,163) -- (250.6,163) .. controls (247.51,163) and (245,160.49) .. (245,157.4) -- cycle ;
\draw    (293,149) -- (314,149) ;
\draw [shift={(316,149)}, rotate = 180] [color={rgb, 255:red, 0; green, 0; blue, 0 }  ][line width=0.75]    (10.93,-3.29) .. controls (6.95,-1.4) and (3.31,-0.3) .. (0,0) .. controls (3.31,0.3) and (6.95,1.4) .. (10.93,3.29)   ;
\draw  [fill={rgb, 255:red, 239; green, 231; blue, 231 }  ,fill opacity=1 ] (244,207.6) .. controls (244,204.51) and (246.51,202) .. (249.6,202) -- (286.4,202) .. controls (289.49,202) and (292,204.51) .. (292,207.6) -- (292,224.4) .. controls (292,227.49) and (289.49,230) .. (286.4,230) -- (249.6,230) .. controls (246.51,230) and (244,227.49) .. (244,224.4) -- cycle ;

\draw (252.6,59.4) node [anchor=north west][inner sep=0.75pt]    {$~~X^{0} \ $};
\draw (244.8,99.4) node [anchor=north west][inner sep=0.75pt]    {$u_{1}^{[ k]} X^{1} \ $};
\draw (123,122.4) node [anchor=north west][inner sep=0.75pt]    {$Y^{[ k]}$};
\draw (393,85.4) node [anchor=north west][inner sep=0.75pt]    {$Y^{[ k+1]} =\ Y^{k} +h\sum _{i=1}^{l} u_{i} X^{i}{\left( Y^{[k]}\right)}$};
\draw (245.8,138.4) node [anchor=north west][inner sep=0.75pt]    {$u_{2}^{[ k]} X^{2} \ $};
\draw (244.8,205.4) node [anchor=north west][inner sep=0.75pt]    {$\ $};
\draw (242.8,204.4) node [anchor=north west][inner sep=0.75pt]    {$u_{m}^{[ k]} X^{m} \ $};

\end{tikzpicture}
}
			\caption{A layer in the neural network control-affine case. Note that the layout is dependent on the choice of temporal discretization. In this case, forward Euler is used, resulting in a ResNet-like structure.}
         	\label{fig:here}
		\end{figure}

	In the sub-Riemannian landmark matching case, we do not consider a weight and bias curve or a neural architecture given by the activation function, but rather a general curve taking values in $\cU$.
 Note, for instance, that if we consider the control-affine setting of \autoref{eq:linaff}, $\cU$ will be equal to $\R^l$, and $v_t = X^0 + \sum_{i=j}^l u_j X^j$. Therefore, 
	\begin{align*}
	    \dot y = X^0 y + \sum_{j=1}^l u_j X^j y. 
	\end{align*}
	Applying the forward Euler method, we obtain 
	\begin{align*}
	    Y_{k}^{[i+1]}=Y_k^{[i]}+h \left( X^0 Y_k^{[i]} +\sum_{j=1}^l u_j^{[i]} X^j Y_k^{[i]}\right),
	\end{align*}
	which can be interpreted as a neural architecture. 
    Note that while in one layer, the nonlinearity is on the data and not on the control variable, 
    the final network is built by composing several layers, using the output from the previous layer as input to the next. 
    Therefore, a nonlinearity is applied to all control variables except to those in the final layer. 
    Note that \citet{Scagliotti2023} has considered a similar neural architecture and proven certain approximation properties in the Euclidean case. 
	A schematic of a layer in the neural network interpretation of sub-Riemannian landmark matching in the setting of control affine systems is seen in \autoref{fig:here}. 
    This type of network has also been studied from a control-theoretic point of view by \citet{Agrachev2021}.

There are, however, some differences between ResNets and sub-Riemannian landmark matching.
Landmarks are typically points in low-dimensional spaces.
But the same framework still applies to points in high-dimensional spaces, such as image data. 
Just as landmarks discretize an image, input data to a ResNet discretizes a larger object. 
This object is  a ``meta-image'', consisting of all possible images of the same class. 
The existence of such an object is the core of the classification problem; shape analysis provides a novel viewpoint.
 We stress that while this interpretation of images is not standard, it is closely related to the notion of \emph{ideas}, discussed in \citet{Ideashape}, where a Gaussian process-based generalization of ResNets was shown to converge to a registration problem of mappings from input and output spaces into spaces of unrealized ``forms'' or ideas, lending its name from Plato's theory of forms.

Another difference between ResNets and landmark matching is how the optimization is performed.
In the deep learning setting, forward and back-propagation are often used \citep[Chapter 5]{bishop:2006:PRML}.
In sub-Riemannian landmark matching, one can either use a shooting method or a gradient flow-type approach, as described in \autoref{sec:compute}. 
The gradient flow algorithm is reminiscent of how a neural network is trained. 
Furthermore, the shooting-based method is also similar to novel particle-based approaches to neural networks, which results in a shooting formulation of ResNets \citep{Vialard2020}. 

Note also that ResNets and sub-Riemannian landmark matching differ in the choice of regularization. 
In the former case, the regularizer is often some matrix norm of the weight matrices \citep[Chapter 5.5]{bishop:2006:PRML}, whereas in the sub-Riemannian landmark matching case the regularization is not directly on the controls, but on the vector fields they parameterize.

In summary, our observations suggest a dictionary between sub-Riemannian landmark matching and ResNets as in \autoref{tbl:resnetdict}.  

\begin{table}[h]
\centering
\begin{tabular}{@{}l l@{}}
\toprule
\textbf{Deep Learning} & \textbf{Landmark matching}          \\ \midrule 
Images                 & Landmarks                           \\ 
Meta-images            & Images                              \\ 
Training network       & Shooting method                     \\ 
Testing                & Warping new landmarks               \\ 
Classification layer   & Forward model                       \\ 
Weights and biases     & Control parameters in $\mathcal{U}$ \\ 
Neural architecture    & ODE Discretization                  \\ \bottomrule
\end{tabular}
\caption{Dictionary between sub-Riemannian landmark matching and ResNets.}\label{tbl:resnetdict}
\end{table}

We remark that geometric control theory can be applied to ResNets just as it can be applied to sub-Riemannian landmark matching. ResNets work because of the nonlinearities. 
If there were none, the generated algebra $\operatorname{Lie}(\cF)$ would be small. 
In fact, it would not contain enough fields to connect two images. 
The nonlinearities assert that the generated Lie algebra is large.
 Further, as observed in the examples, the inclusion of meaningful nonlinearities means that the decreased time step, directly translating to an increase in the number of layers, is worthwhile.

\section{Conclusions}\label{sec:conclusions}

There are two main ideas in this paper. 
The first is that sub-Riemannian landmark matching can be approached from a geometric point of view, and that it is possible to derive equations governing the evolution of the control parameters directly from variational principles.
In \autoref{sec:proof} we also prove the existence of minimizers.

The second idea concerns the connection between shape analysis and deep learning. 
Shape analysis offers a view of deep learning as a way to find parameterized warps for moving initial landmarks to target landmarks, a problem with inherent geometry. 
In particular, it explains the importance of the nonlinearities (activation functions) of a neural network: these ensure that the bracket generated algebra is large so that the range of transformations from input to output is large.
An interesting opening here is to investigate more closely how the choice of activation functions affect the bracket generated algebra; to find nonlinearities that maximize it is a plausible network design objective.

The fact that geometric and control-theoretic methods can be applied to understand deep learning methods seem to imply that there is more to be done. 
A future direction of research could be to further the understanding of deep neural networks as geometric objects in their own right.

\appendix
\section{Geometric preliminaries}

\label{sec:geometry}
In this section, we provide a brief and informal introduction to concepts from differential geometry used in the paper. 
The reader already familiar with this field can directly skip to \autoref{sec:landmark}. 
For a more substantial introduction to geometry, see for instance \citet{Lee1997}.

We begin with an intuitive idea of what a manifold is. 
A manifold $M$ of dimension $n$ is a space which locally resembles Euclidean space $\R^n$.
One can picture this as a generalization of an embedded surface. 
The set of vectors tangent to the surface at a point $p \in M$ is a vector space called the \emph{tangent space at $p$} and is denoted $T_p M$, and the disjoint union of the tangent spaces at each point of the manifold is called the \emph{tangent bundle} of $M$ and is denoted by $TM$. 
A vector field $v$ on $M$ is an assignment of a vector $v(p)\in T_pM$ to each point $p\in M$.
The space of all smooth vector fields is denoted $\mathfrak{X}(M)$. 

A \emph{Riemannian manifold} is a manifold equipped with a \emph{Riemannian metric} $g$, i.e., at each point $p\in M$ the metric $g_p$ is an inner product on $T_p M $. A Riemannian metric allows us to define distances, angles, and curvature of the manifold. 

A \emph{Lie group} $G$ is a manifold as well as a multiplicative group with the property that the map from $G\times G$ to $G$ given by $(g,h) \mapsto gh^{-1}$ is smooth. 
The tangent space of $G$ at the identity $e\in G$ is referred to as the \emph{Lie algebra} $\mathfrak{g}$ of $G$.
It is equipped with a \emph{Lie bracket} mapping $\mathfrak{g} \times \mathfrak{g} \to \mathfrak{g}$, denoted by $[v,w]$.
A linear subspace closed under the Lie bracket is called a \emph{Lie subalgebra}.

A \emph{distribution} $\cD$ of $M$ is a collection of linear subspaces $\cD_p \subset T_p M$ at each point $p\in M$. 
If it is closed under the Lie bracket, then the distribution is said to be \emph{integrable}, and it defines, at least locally, a \emph{foliation} of $M$. 
It  can be thought of a partition of $M$ into submanifolds called the \emph{leaves} of the foliation.  

Having briefly mentioned the geometric tools needed to read the paper, we can move on to introduce concepts from shape analysis. 
\section{Existence of minimizers}
\label{sec:proof}
In this appendix, we first formalize the analytical setting, for which we then prove the existence of solutions to the sub-Riemannian landmark matching problem.

Let $\cU$ and $\cH$ be Hilbert spaces. Assume that $\cU$ is separable and that $\cH$ is 
continuously embedded into $C^1(M;TM)$. 
Examples of admissible spaces  are the Sobolev spaces of high enough regularity that were used in \autoref{constrained_landmarks}. 

The variables parameterizing the vector fields will be in the space
\begin{align*}
    \mathsf{H}^k([0,1],\cU),
\end{align*}
i.e. the space of all functions in $\mathsf{L}^2([0,1],\cU)$ such that the weak derivatives up to $k\geq 1$ also are in $\mathsf{L}^2([0,1],\cU)$.  
 Here $\mathsf{L}^2([0,1],\cU)$ denotes the space of Bochner measurable square integrable functions taking values in $\cU$. 
Note that 
\begin{align*}
   \mathsf{H}^0([0,1],\cU)=\mathsf{L}^2([0,1],\cU). 
\end{align*}
See \citet{veraar} for more details on Bochner spaces and Banach space-valued Sobolev spaces. 

The vector fields are in $\mathsf{L}^2([0,1],\cH)$ and will be parameterized by $u \in \mathsf{H}^k([0,1],\cU)$ via a function  
\begin{equation*}
    F\colon\mathsf{H}^k([0,1],\cU) \to \mathsf{L}^2([0,1],\cH).
\end{equation*}
We assume that $F$ is weakly continuous. 

The set of parametrized vector fields formally determines a subset $\cS \subset \mathfrak{X}(M)$ given by \autoref{eq:S}. 
In general, $\cS$ is not a Lie subalgebra since the Lie bracket is not closed on $\cS$.

Since $\dot{\gamma}(t) = v(t) \circ \gamma(t)$ and $v(t) = F(u)_t$ the equation governing the trajectory $y(t)$ of points is  given by 
\begin{align}
    y(t) =y(0)+\int_0^t F(u)_s(y(s)) \dd s.
\end{align}
The integral $\int_0^t F(u)_s(y(s)) \dd s$ is bounded since $\cH$ 
is continuously embedded into $C^1(M;TM)$.

The minimization problem is therefore 
\begin{align}
    	&\!\min_{u \in \mathsf{H}^k([0,1];\cH)} \quad \sum_{i=1}^m d_N(h(y_i)),c_i)^2+\mathcal{R}(u)\label{eq:targetvlm}\\
			&\text{s.t.  } \quad y_i(t) =x _i +\int_0^t F(u)_s(y_i(s)) \dd s  \label{eq:constvlm}.
\end{align}
where 
\begin{align*}
    \mathcal{R}({u})=\int_0^1 \mathpzc{R}(u_t) \mathrm{d}t.
\end{align*} 
Here $\mathpzc{R}$ is a map from $\cU$ into $\R$. 
Existence of minimizers to problem \eqref{eq:targetvlm}-\eqref{eq:constvlm} is given by the following theorem:
\begin{theorem}
\label{th:exst}
Assume that 
$\mathcal{R}(u)$ is lower semi-continuous and that $\mathpzc{R}(u_t)$ bounds the $\cU$-norm for all weak derivatives up to order $k$.
If the forward model $h\colon M\to N$ is continuous,  then there exists a minimizer $\tilde{u}\in \mathsf{H}^k([0,1],\cU)$ to problem \eqref{eq:targetvlm}-\eqref{eq:constvlm}. 
\end{theorem}
To prove this theorem, we shall use the direct method of variational calculus \citep{Bruveris,thesisglaunes}. 
\begin{proof}
Note first that $y_i(1)= \gamma_1(x_i)$, the flow of the time-dependent vector field $F(u)$ at the time $1$ evaluated at $x_i$. 
We let $u^n \in \mathsf{H}^k([0,1],\cU)$ be a minimizing sequence, that is to say,  
\begin{align*}
E(u^n) \to \inf_u E(u).
\end{align*}
By assumption there is a constant $c>0$ such that  $\|u^n\|_{\mathsf{H}^k([0,1],\cU)}^2 \leq c \mathcal{R}(u) $. 
It holds that 
\begin{align*}
\|u^n\|_{\mathsf{H}^k([0,1],\cU)}^2 \leq c\mathcal{R}(u^n)\leq C E(u^n )
\end{align*}     
and so $u^n$ is a bounded sequence in $\mathsf{H}^k([0,1],\cU)$ and we can extract a subsequence $(u^{n_k})$ converging weakly to $\tilde{u}$. 
Note that $F$ will define a sequence of vector fields $v^{n_k}=F(u^{n_k})$ in  $\mathsf{L}^2([0,1],\cH)$.  The sequence $v^{n_k}$ will converge weakly to $\tilde{v} = F(\tilde{u})$ since $F$ is weakly continuous. It is now possible to proceed as in \citet[Lemma 5]{thesisglaunes} to obtain pointwise convergence of the flows. 
By the assumption that the regularizer $\mathcal{R}(u)$ is lower semi-continuous, and that the forward model is continuous, it follows that 
\begin{align*}
&\inf_{u \in \mathsf{L}^2([0,1]\cH)} E(u) \leq E(\tilde{u})=\mathcal{R}(\tilde{u})+ \sum_{i=1}^m d_N(h(\tilde\gamma_1(x_i)),c_i)^2\\
&\leq \liminf_{n_k \to \infty }\mathcal{R}(u^{n_k})+\lim_{n_k \to \infty } \sum_{i=1}^m d_N(h(\tilde\phi_1^{n_k}(x_i)),c_i)^2\leq \inf_{u \in \mathsf{L}^2([0,1]\cH)} E(u).
\end{align*}
This concludes the proof.
\end{proof}

\bibliographystyle{amsplainnat.bst}
\bibliography{shprnt_refs.bib}

\begin{thebibliography}{45}
\providecommand{\natexlab}[1]{#1}
\providecommand{\url}[1]{\texttt{#1}}
\providecommand{\urlprefix}{URL }
\providecommand{\eprint}[2][]{\url{#2}}

\bibitem[{Agrachev and Sachkov(2004)}]{Agrachev2004}
A.~A. Agrachev and Y.~L. Sachkov, \emph{Control Theory from the Geometric
  Viewpoint}, Springer, 2004.

\bibitem[{Agrachev and Sarychev(2021)}]{Agrachev2021}
A.~A. Agrachev and A.~Sarychev, {Control on the Manifolds of Mappings with a
  View to the Deep Learning}, \emph{J. Dyn. Control Syst.}  (2021).

\bibitem[{Arguill{è}re et~al.(2015)Arguill{è}re, Tr{\'e}lat, Trouv{\'e}, and
  Younes}]{arguillere2015}
S.~Arguill{è}re, E.~Tr{\'e}lat, A.~Trouv{\'e}, and L.~Younes, {Shape
  deformation analysis from the optimal control viewpoint}, \emph{J. Math.
  Pures. Appl.} \textbf{104} (2015), 139--178.

\bibitem[{Arguill{è}re et~al.(2016)Arguill{è}re, Trélat, Trouvé, and
  Younes}]{Sylvain2016}
S.~Arguill{è}re, E.~Trélat, A.~Trouvé, and L.~Younes, Registration of
  multiple shapes using constrained optimal control, \emph{SIAM J. Imaging
  Sci.} \textbf{9} (2016), 344--385.

\bibitem[{Beauchard et~al.(2023)Beauchard, Le~Borgne, and
  Marbach}]{Beauchard2023}
K.~Beauchard, J.~Le~Borgne, and F.~Marbach, On expansions for nonlinear systems
  error estimates and convergence issues, \emph{Comptes Rendus. Mathématique}
  \textbf{361} (2023), 97–189.

\bibitem[{Beg et~al.(2005)Beg, Miller, Trouv{\'{e}}, and Younes}]{Beg2005}
M.~F. Beg, M.~I. Miller, A.~Trouv{\'{e}}, and L.~Younes, Computing large
  deformation metric mappings via geodesic flows of diffeomorphisms, \emph{Int.
  J. Comput. Vis.} \textbf{61} (2005), 139--157.

\bibitem[{Ben~Amor et~al.(2023)Ben~Amor, Arguillere, and
  Shao}]{DBLP:journals/corr/abs-2102-07951}
B.~Ben~Amor, S.~Arguillere, and L.~Shao, {ResNet-LDDMM}: Advancing the {LDDMM
  }framework using deep residual networks, \emph{IEEE Trans. Pattern Anal.
  Mach. Intell.} \textbf{45} (2023), 3707--3720.

\bibitem[{Bishop(2006)}]{bishop:2006:PRML}
C.~M. Bishop, \emph{Pattern Recognition and Machine Learning}, Springer, 2006.

\bibitem[{Bistoquet et~al.(2008)Bistoquet, Oshinski, and Skrinjar}]{heartmri}
A.~Bistoquet, J.~Oshinski, and O.~Skrinjar, Myocardial deformation recovery
  from cine {MRI} using a nearly incompressible biventricular model, \emph{Med.
  Image Anal.} \textbf{12} (2008), 69--85.

\bibitem[{Bruveris(2012)}]{Bruveris}
M.~Bruveris, \emph{Geometry of Diffeomorphism Groups and Shape Matching}, Ph.D.
  thesis, Imperial College London, 2012.

\bibitem[{Bruveris and Holm(2015)}]{Bruveris2013}
M.~Bruveris and D.~D. Holm, \emph{Geometry of Image Registration: The
  Diffeomorphism Group and Momentum Maps}, p. 19–56, Springer New York, 2015.

\bibitem[{Celledoni et~al.(2021)Celledoni, Ehrhardt, Etmann, Owren, Schönlieb,
  and Sherry}]{Celledoni2020}
E.~Celledoni, M.~J. Ehrhardt, R.~I. Etmann, Cand~McLachlan, B.~Owren, C.-B.
  Schönlieb, and F.~Sherry, Structure-preserving deep learning, \emph{Eur. J.
  Appl. Math.} \textbf{32} (2021), 888–936.

\bibitem[{Celledoni et~al.(2023)Celledoni, Gl\"{o}ckner, Riseth, and
  Schmeding}]{Celledoni2023}
E.~Celledoni, H.~Gl\"{o}ckner, J.~N. Riseth, and A.~Schmeding, Deep neural
  networks on diffeomorphism groups for optimal shape reparametrization,
  \emph{{BIT} Numer. Math.} \textbf{63} (2023).

\bibitem[{Ceritoglu et~al.(2013)Ceritoglu, Tang, Chow, Hadjiabadi, Shah, Brown,
  Burhanullah, Trinh, Hsu, Ament, Crocetti, Mori, Mostofsky, Yantis, Miller,
  and Ratnanather}]{brainC}
C.~Ceritoglu, X.~Tang, M.~Chow, D.~Hadjiabadi, D.~Shah, T.~Brown,
  M.~Burhanullah, H.~Trinh, J.~Hsu, K.~Ament, D.~Crocetti, S.~Mori,
  S.~Mostofsky, S.~Yantis, M.~I. Miller, and T.~J. Ratnanather, Computational
  analysis of {LDDMM} for brain mapping, \emph{Front. Neurosci.} \textbf{7}
  (2013), 151.

\bibitem[{Cifor et~al.(2013)Cifor, Risser, Chung, Anderson, and
  Schnabel}]{cancer}
A.~Cifor, L.~Risser, D.~Chung, E.~M. Anderson, and J.~A. Schnabel, Hybrid
  feature-based diffeomorphic registration for tumor tracking in {2-D} liver
  ultrasound images, \emph{IEEE Trans. Med. Imaging} \textbf{32} (2013),
  1647--1656.

\bibitem[{Ganaba(2021)}]{nader2021}
N.~Ganaba, Deep learning: {H}ydrodynamics, and {L}ie-{P}oisson
  {H}amilton-{J}acobi theory, arXiv:2105.09542v2, 2021.

\bibitem[{Glaunes(2005)}]{thesisglaunes}
J.~A. Glaunes, \emph{Transport par difféomorphismes de points, demesures et de
  courants pour la comparaisonde formes et l’anatomie numérique.}, Ph.D.
  thesis, Université Paris 13, 2005.

\bibitem[{Gris et~al.(2018)Gris, Durrleman, and Trouv{\'e}}]{gris1}
B.~Gris, S.~Durrleman, and A.~Trouv{\'e}, A {Sub-Riemannian} modular framework
  for diffeomorphism-based analysis of shape ensembles, \emph{SIAM J. Imaging
  Sci.} \textbf{11} (2018), 802--833.

\bibitem[{Hamilton(1982)}]{Hamilton1982}
R.~S. Hamilton, The inverse function theorem of nash and moser, \emph{Bull. Am.
  Math. Soc.} \textbf{7} (1982), 65–222.

\bibitem[{He et~al.(2016)He, Zhang, Ren, and Sun}]{He2015}
K.~He, X.~Zhang, S.~Ren, and J.~Sun, Deep residual learning for image
  recognition, \emph{2016 IEEE Conference on Computer Vision and Pattern
  Recognition (CVPR)}, pp. 770--778, 2016.

\bibitem[{H\"{o}rmander(1967)}]{Hrmander1967}
L.~H\"{o}rmander, Hypoelliptic second order differential equations, \emph{Acta
  Math.} \textbf{119} (1967), 147--171.

\bibitem[{Hyt{\"o}nen et~al.(2016)Hyt{\"o}nen, van Neerven, Veraar, and
  Weis}]{veraar}
T.~Hyt{\"o}nen, J.~van Neerven, M.~Veraar, and L.~Weis, \emph{Analysis in
  Banach Spaces: Volume I: Martingales and Littlewood-Paley Theory}, Springer,
  2016.

\bibitem[{Joshi and Miller(2000)}]{Joshi2000}
S.~C. Joshi and M.~I. Miller, Landmark matching via large deformation
  diffeomorphisms, \emph{{IEEE} Trans. Image Process.} \textbf{9} (2000),
  1357--1370.

\bibitem[{Jurdjevic(1999)}]{Jurdjevic1999}
V.~Jurdjevic, Optimal control, geometry, and mechanics, \emph{Mathematical
  Control Theory}, pp. 227--267, Springer, 1999.

\bibitem[{Lee(1997)}]{Lee1997}
J.~M. Lee, \emph{Riemannian Manifolds}, Springer, 1997.

\bibitem[{Li et~al.(2017)Li, Chen, Cheng, and Weinan}]{Li2017}
Q.~Li, L.~Chen, T.~Cheng, and E.~Weinan, {Maximum Principle Based Algorithms
  for Deep Learning}, \emph{J. Mach. Learn. Res.} \textbf{18} (2017),
  5998–6026.

\bibitem[{Mansi et~al.(2010)Mansi, Pennec, Sermesant, Delingette, and
  Ayache}]{heartmri2}
T.~Mansi, X.~Pennec, M.~Sermesant, H.~Delingette, and N.~Ayache, {iLogDemons}:
  A demons-based registration algorithm for~tracking incompressible elastic
  biological tissues, \emph{Int. J. Comput. Vis.} \textbf{92} (2010), 92--111.

\bibitem[{Marsden and Ratiu(1999)}]{Marsden1999}
J.~E. Marsden and T.~Ratiu, \emph{Introduction to Mechanics and Symmetry},
  Springer, 1999.

\bibitem[{Miller(1974)}]{Webb}
W.~Miller, Computational complexity and numerical stability, \emph{{Proceedings
  of the Sixth Annual ACM Symposium on Theory of Computing}}, p. 317–322,
  Association for Computing Machinery, 1974.

\bibitem[{Mogensen and Riseth(2018)}]{Optim.jl-2018}
P.~K. Mogensen and A.~N. Riseth, Optim: A mathematical optimization package for
  {Julia}, \emph{J. Open Source Softw.} \textbf{3} (2018), 615.

\bibitem[{Nielsen~Clelland et~al.(2009)Nielsen~Clelland, Moseley, and
  Wilkens}]{Clelland2009}
J.~Nielsen~Clelland, C.~Moseley, and G.~R. Wilkens, Geometry of control-affine
  systems, \emph{Symmetry Integr. Geom.} \textbf{5} (2009), 095.

\bibitem[{Owhadi(2023)}]{Ideashape}
H.~Owhadi, Do ideas have shape? idea registration as the continuous limit of
  artificial neural networks, \emph{Physica D} \textbf{444} (2023), 133592.

\bibitem[{Paszke et~al.(2019)Paszke, Gross, Massa, Lerer, Bradbury, Chanan,
  Killeen, Lin, Gimelshein, Antiga, Desmaison, Kopf, Yang, DeVito, Raison,
  Tejani, Chilamkurthy, Steiner, Fang, Bai, and Chintala}]{NEURIPS2019_9015}
A.~Paszke, S.~Gross, F.~Massa, A.~Lerer, J.~Bradbury, G.~Chanan, T.~Killeen,
  Z.~Lin, N.~Gimelshein, L.~Antiga, A.~Desmaison, A.~Kopf, E.~Yang, Z.~DeVito,
  M.~Raison, A.~Tejani, S.~Chilamkurthy, B.~Steiner, L.~Fang, J.~Bai, and
  S.~Chintala, {PyTorch}: {An Imperative Style, High-Performance Deep Learning
  Library}, \emph{Advances in Neural Information Processing Systems 32}, pp.
  8024--8035, Curran Associates, Inc., 2019.

\bibitem[{Rackauckas and Nie(2017)}]{rackauckas2017differentialequations}
C.~Rackauckas and Q.~Nie, {Differentialequations.jl}--a performant and
  feature-rich ecosystem for solving differential equations in {J}ulia,
  \emph{J. Open Res. Softw.} \textbf{5} (2017).

\bibitem[{{Revels} et~al.(2016){Revels}, {Lubin}, and
  {Papamarkou}}]{RevelsLubinPapamarkou2016}
J.~{Revels}, M.~{Lubin}, and T.~{Papamarkou}, Forward-mode automatic
  differentiation in {J}ulia, \emph{CoRR}  (2016).

\bibitem[{Risser et~al.(2013)Risser, Vialard, Baluwala, and Schnabel}]{lung}
L.~Risser, F.~X. Vialard, H.~Y. Baluwala, and J.~A. Schnabel,
  Piecewise-diffeomorphic image registration: {A}pplication to the motion
  estimation between {3D} {CT} lung images with sliding conditions, \emph{Med.
  Image Anal.} \textbf{17} (2013), 182--193.

\bibitem[{Scagliotti(2023)}]{Scagliotti2023}
A.~Scagliotti, Deep learning approximation of diffeomorphisms via
  linear-control systems, \emph{Math. Control. Relat. Fields} \textbf{13}
  (2023), 1226--1257.

\bibitem[{Sussmann(1986)}]{Sussmann1986}
H.~J. Sussmann, A product expansion for the {C}hen series, C.~I. Byrnes and
  A.~Lindquist, eds., \emph{{Theory and Applications of Nonlinear Control
  Systems}}, p. 323–335, North-Holland, 1986.

\bibitem[{Tabuada and Gharesifard(2021)}]{Tabuada2021UniversalAP}
P.~Tabuada and B.~Gharesifard, Universal approximation power of deep residual
  neural networks via nonlinear control theory, \emph{ICLR}, 2021.

\bibitem[{Thompson(1992)}]{Thompson1992}
D.~Thompson, \emph{On Growth and Form}, Cambridge University Press, 1992.

\bibitem[{Vialard et~al.(2020)Vialard, Kwitt, Wei, and
  Niethammer}]{Vialard2020}
F.~Vialard, R.~Kwitt, S.~Wei, and M.~Niethammer, A shooting formulation of deep
  learning, H.~Larochelle, M.~Ranzato, R.~Hadsell, M.~F. Balcan, and H.~Lin,
  eds., \emph{Advances in Neural Information Processing Systems}, vol.~33, pp.
  11828--11838, Curran Associates, Inc., 2020.

\bibitem[{{Younes}(2010)}]{Younes2010}
L.~{Younes}, \emph{{S}hapes and {D}iffeomorphisms}, Springer, 2010.

\bibitem[{Younes(2020)}]{Younes2020b}
L.~Younes, Diffeomorphic learning, \emph{J. Mach. Learn. Res.} \textbf{21}
  (2020), 1--28.

\bibitem[{Younes et~al.(2020)Younes, Gris, and Trouv{\'e}}]{Younes2020}
L.~Younes, B.~Gris, and A.~Trouv{\'e}, \emph{Sub-Riemannian Methods in Shape
  Analysis}, pp. 463--495, Springer, 2020.

\bibitem[{{Ö}ktem et~al.(2016){Ö}ktem, Chong, Nevzat, Pradeep, and
  Chandrajit}]{oktem2016}
O.~{Ö}ktem, C.~Chong, D.~Nevzat, R.~Pradeep, and B.~Chandrajit, Shape based
  image reconstruction using linearized deformations, \emph{Inverse Probl.}
  \textbf{33} (2016).

\end{thebibliography}

\end{document}